\documentclass{conm-p-l} 
\usepackage{amssymb,array,delarray,latexsym,epsfig} 

\newlength{\sh}
\newlength{\baker}
\newlength{\greg}
\newlength{\fw}
\newlength{\jmr}
\newlength{\jfc}
\newlength{\bernd}
\newlength{\jones}
\newlength{\mati}
\newlength{\tung}
\newlength{\sil}
\newlength{\koi}
\newlength{\gala}
\newtheorem{kho}{Khovanski's Theorem on Real Fewnomials} 

\newtheorem{jst}{The JST Theorem}

\newtheorem{tungy}{Tung's Theorem}

\newtheorem{lemma}{Lemma}

\newtheorem{dfn}{Definition}

\newtheorem{main}{theorem} 
\newtheorem{thm}[main]{Theorem}
\newtheorem{cor}{Corollary}
\newtheorem{rem}{Remark}	
\newtheorem{ex}{Example}

\newcommand{\tf}{{\tilde{f}}}
\newcommand{\twF}{{\tilde{F}}} 
 
\renewcommand{\mod}{\mathbf{mod}} 
\newcommand{\pspa}{{\mathbf{PSPACE}}}

\newcommand{\am}{{\mathbf{AM}}} 

\newcommand{\np}{{\mathbf{NP}}}
\newcommand{\conp}{{\mathbf{coNP}}}

\newcommand{\bpp}{{\mathbf{BPP}}}
\newcommand{\crap}{\pp^{\np^\np}}
\newcommand{\pp}{\mathbf{P}}
\newcommand{\hn}{\mathbf{HN}}
\newcommand{\nc}{\mathbf{NC}}
\newcommand{\expt}{{\mathbf{EXPTIME}}}
\newcommand{\eps}{\varepsilon}
\newcommand{\cA}{\mathcal{A}}

\newcommand{\cO}{\mathcal{O}}
\newcommand{\supp}{\mathrm{Supp}}
\newcommand{\conv}{\mathrm{Conv}}
\newcommand{\size}{\mathrm{size}}
\newcommand{\thth}{{\underline{\mathrm{th}}}}
\newcommand{\rd}{ {\underline{ \mathrm{rd} } }  }
\newcommand{\st}{ {\underline{ \mathrm{st} } }  }
\newcommand{\nd}{{\underline{\mathrm{nd}}}}
\newcommand{\Pro}{{\mathbb{P}}}

\newcommand{\Q}{\mathbb{Q}}
\newcommand{\R}{\mathbb{R}}
\newcommand{\C}{\mathbb{C}}
\newcommand{\N}{\mathbb{N}}
\newcommand{\Z}{\mathbb{Z}}

\newcommand{\pert}{\mathrm{Pert}}

\newcommand{\res}{\mathrm{Res}}

\newcommand{\Zn}{\Z^n}

\newcommand{\Rn}{\R^n}

\newcommand{\Cn}{\C^n}

\newcommand{\Csn}{{(\C^*)}^n}
\renewcommand{\qed}{$\blacksquare$}
\newcommand{\cM}{{\mathcal{M}}}

\newcommand{\cI}{\mathcal{I}}
\newcommand{\cR}{\mathcal{R}}
\newcommand{\cP}{\mathcal{P}} 
\newcommand{\cQ}{\mathcal{Q}}
\newcommand{\cS}{\mathcal{S}}
\newcommand{\cC}{\mathcal{C}}
\newcommand{\cU}{\mathcal{U}}
\newcommand{\bO}{\mathbf{O}}
\newcommand{\vol}{\mathrm{Vol}}
\newcommand{\htp}{\mathrm{HTP}}
\newcommand{\ratcurve}{\mathrm{RatCurve}}
\newcommand{\biggy}{\mathrm{Big}}

\begin{document}

\title[Algebraic Geometry Over Four Rings]{ 
Algebraic Geometry Over Four Rings and the Frontier to Tractability}  

\author{J.\ Maurice Rojas}\thanks{To appear in a volume of Contemporary 
Mathematics:   
Proceedings of a Conference on Hilbert's Tenth Problem and Related Subjects 
(University of Gent, November 1--5, 1999), edited by 
Jan Denef, Leonard Lipschitz, Thanases Pheidas, and Jan Van 
Geel, AMS Press. 
This research was partially supported by a 
Hong Kong CERG Grant.}   

\address{Department of Mathematics, City University of Hong Kong, 
83 Tat Chee Avenue, Kowloon, HONG KONG} 
\email{mamrojas@math.cityu.edu.hk\\ {\it Web-Page:}  
http://www.cityu.edu.hk/ma/staff/rojas } 

\dedicatory{This paper is dedicated to Steve Smale on the 
occasion of his $70^\thth$ birthday.}  

\date{\today} 

\begin{abstract} 
We present some new and recent algorithmic results concerning polynomial 
system solving over various rings. In particular, we 
present some of the best recent bounds on:  
\begin{itemize}
\item[(a)]{\mbox{the complexity of calculating the complex dimension of an 
algebraic set} } 
\item[(b)]{the height of the zero-dimensional part of an algebraic set over 
$\C$} 
\item[(c)]{the number of connected components of a semi-algebraic set}
\end{itemize} 
We also present some results which significantly lower the complexity of 
deciding the emptiness of hypersurface intersections over $\C$ and $\Q$, 
given the truth of the Generalized Riemann Hypothesis. 
Furthermore, we state some recent progress on the decidability of the 
prefixes $\exists\forall\exists$ and $\exists\exists\forall\exists$, 
quantified over the positive integers. As an application, we conclude with 
a result connecting Hilbert's Tenth Problem in three variables and height 
bounds for integral points on algebraic curves. 

This paper is based on three lectures presented at the 
conference corresponding to this proceedings volume. 
The titles of the lectures were ``Some Speed-Ups in Computational 
Algebraic Geometry,'' ``Diophantine Problems Nearly in the Polynomial 
Hierarchy,'' and ``Curves, Surfaces, and the Frontier to Undecidability.''
\end{abstract} 

\mbox{}\\
\vspace{-.3in}
\maketitle

\mbox{}\\
\vspace{-.7in}
\tableofcontents 

\mbox{}\\
\vspace{-.7in}
\section{Introduction}
\label{sec:intro} 
This paper presents an assortment of algorithmic and combinatorial 
results that the author hopes is useful to 
experts in arithmetic geometry and diophantine complexity. 
While the selection of results may appear somewhat eclectic, 
there is an underlying motivation: determining the boundary 
to tractability for polynomial equation solving in various settings. 
The notion of tractability here will mean membership in a particular 
well-known complexity class depending on the underlying ring and input 
encoding. As an example of this principle, we point out that our brief tour 
culminates with a result giving evidence for the following assertion: The 
recursive unsolvability of deciding the existence of integral roots 
for multivariate polynomials begins with polynomials in {\bf three} variables. 
The sharpest current threshold is still nine variables (for 
{\bf positive} integral roots) \cite{jones9}.\footnote{James P.\ Jones, 
the author of \cite{jones9}, attributes the nine variables result to 
Matiyasevich.} 

Our main results will first be separated into the underlying 
ring of interest, here either $\C$, $\R$, $\Q$, or $\Z$. 
Within each group of results, we will warm up with a 
non-trivial result involving univariate polynomials. All 
necessary proofs are elaborated in section \ref{sec:proofs}, and 
our main underlying computational models will either be the 
classical {\bf Turing machine} \cite{papa} or the 
{\bf BSS machine over $\pmb{\C}$} \cite{bcss}. 
The two aforementioned references are excellent sources for 
further complexity-theoretic background, but we will only require a minimal 
acquaintance with these computational models. 

Before embarking on the full technical statements 
of our main theorems, let us see some concrete examples 
to whet the readers appetite, and further ground the 
definitions we will later require. 

\subsection{A Sparse $\pmb{3\times 3}$ Polynomial System}  
\label{sub:3by3}
The solution of sparse polynomial systems is a problem with   
numerous applications outside, as well as inside, mathematics. The 
analysis of chemical reactions \cite{gaterhub} and the computation of 
equilibria in game-theoretic models \cite{mucks} are but two diverse examples. 

More concretely, consider the following system of $3$ polynomial equations in 
$3$ variables: 
\begin{eqnarray}
\label{eq:3by3}
144+2x-3y^2+x^7y^8z^9 &=& 0 \notag \\
-51+5x^2-27z+x^9y^7z^8 &=& 0 \\
7-6x+8x^8y^9z^7-12x^8y^8z^7 &=& 0. \notag 
\end{eqnarray} 
Let us see if the system (\ref{eq:3by3}) has any {\bf complex} 
roots and, if so, count how many there are. Any terminology or results applied 
here will be clarified further in section \ref{sec:complex}. 

Note that the total degree\footnote{ The {\bf total degree} of a polynomial is 
just the maximum of the sum of the exponents in any monomial term of the 
polynomial.} of each polynomial above is 24. By 
an 18$^\thth$-century theorem of \'Etienne B\'ezout \cite{shafa}, 
we can bound from 
above the number of complex roots of (\ref{eq:3by3}), assuming 
this number is finite, by $24\cdot 24 \cdot 24 = \mathbf{13824}$. 
However, a more precise 20$^\thth$-century bound can be obtained 
by paying closer attention to the monomial term structure 
of (\ref{eq:3by3}): Considering the 
convex hull of\footnote{i.e., smallest convex set in $\R^3$ containing...} 
the exponent vectors of each equation in 
(\ref{eq:3by3}), one obtains three tetrahedra. 
\mbox{}\hspace{1cm}\epsfig{file=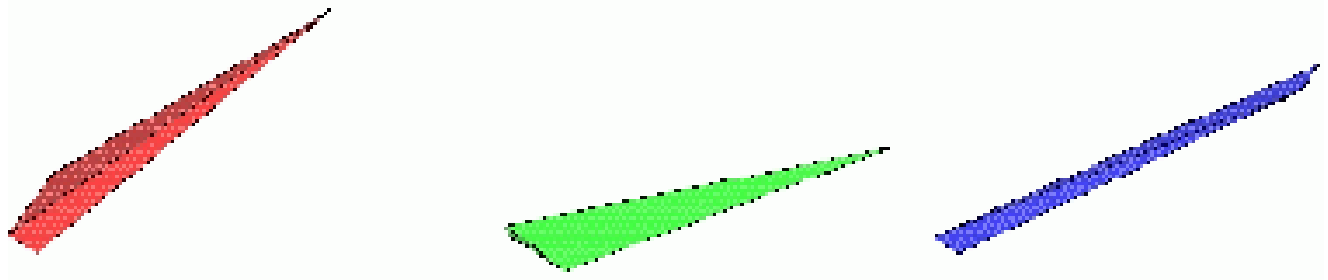,height=1.2in}

These are the {\bf Newton polytopes} of (\ref{eq:3by3}), and 
their {\bf mixed volume}, by a beautiful theorem of 
David N.\ Bernshtein from the 1970's \cite{bernie}, turns out to be a 
much better upper bound on the number of complex roots (assuming there 
are only finitely many). For our polynomial system (1), this bound 
is\footnote{ \label{see} Please see the Appendix for further details on the 
theory and implementation behind our examples. }  {\bf 145}. 

Now to decide whether (\ref{eq:3by3}) has any complex roots, 
we can attempt to find a univariate polynomial whose roots 
are some simple function of the roots of (\ref{eq:3by3}). {\bf Elimination 
theory} allows one to do this, and a particularly effective combinatorial 
algorithm is given in theorem \ref{main:complex} of section 
\ref{sec:complex}. For example, the roots of 
\tiny
\vspace{-.4cm}
\begin{center}
\[ \pmb{P(u)}:= 268435456 u^{145}    -138160373760 u^{137}    
-30953963520 u^{130} +3446308601856 u^{129}    -25165824000 u^{123}  \] 
\vspace{-.4cm}  
\[ -26293995307008 u^{122}  -1694282972921856 u^{121}    +323419618934784 
u^{120} -6995155353600 u^{115} \]  
\vspace{-.4cm}  
\[+87379566133248 u^{114}   +10198949486395392 u^{113}    
-166099501774798848 u^{112}    -112538419200 u^{108}\]
\vspace{-.4cm}  
\[    -82834929745920 u^{107} -324798104395579392 u^{106}
-4419977097552592896 u^{105}   +589824000000 u^{101}\] 
\vspace{-.4cm}  
\[   -35724722176000 u^{100} +8364740005330944 u^{99}   
+4439548695657775104 u^{98} -26917017845238005760 u^{97}\]  
\vspace{-.4cm}  
\[+37910937600000 u^{93}  +51523633570381824 u^{92} 
 -1791672886920019968 u^{91} -848160250027183521792 u^{90} \]
\vspace{-.4cm}  
\[ +616996999355281440768 u^{89}   -664995358310400 u^{85} 
+1524560547831644160 u^{84} +745863497970172674048 u^{83} \]
\vspace{-.4cm}  
\[ +17539603347891497287680 u^{82}+994210006214153207808 u^{81}
+12899450880000 u^{78}   
-47322888233287680 u^{77}\]
\vspace{-.4cm}  
\[+33981667956844904448 u^{76}-4986502987101813633024 u^{75}
+119063825168001672019968 u^{74}\]
\vspace{-.4cm}  
\[+31576057329392164012032 u^{73}+751796121600000 u^{70}
-9866721074229006336 u^{69}\]
\vspace{-.4cm}  
\[+1882463818496535244800 u^{68}+3052871408440654112816640 u^{67}
+380423482789919103664128 u^{66}\]
\vspace{-.4cm}  
\[+34866943014558674976768 u^{65} +279569449114214400 u^{62}
-302173847078728854528u^{61}\]
\vspace{-.4cm}  
\[-534702070464812022223872 u^{60} -14973258769647086979053568 u^{59}
+4994218012036588712165376 u^{58}\]
\vspace{-.4cm}  
\[-2021795433676800 u^{55}+8296585706519424000 u^{54}
+25005465159580886376960 u^{53} -3783799262749190677321536 u^{52}\]
\vspace{-.4cm}  
\[+35916388899232830509942784 u^{51}+6316741393466865886715904 u^{50}
-61674073526016000 u^{47}\]
\vspace{-.4cm}  
\[ -554525302200721744896 u^{46}+812163230435877273319104 u^{45}
-\underline{2947435596503653060289376000} u^{44}\]
\vspace{-.4cm}  
\[-141780781258618244980543488 u^{43}+ 6318299549796897024 u^{39}
-41096279946826872821088 u^{38}\]
\vspace{-.4cm}  
\[+294236770231877581913540688 u^{37}+326253143719924635239730432 u^{36}
-8845750586564412369214464 u^{35}\]
\vspace{-.4cm}  
\[ -29428437386188800 u^{32}+886156671237883112160 u^{31} 
-12033942692990286448093392 u^{30} \cdots \] 
\vspace{-.4cm}  
\[-21345681203414534849440320u^{29}+176061998413186705562222592 u^{28}
-8770384173478164480 u^{24}\]
\vspace{-.4cm}  
\[+258178048486605790963020 u^{23} +482019749452059431164020 u^{22}
-11741024693522572606851840 u^{21}\]
\vspace{-.4cm}  
\[+32803667644608000 u^{17}-3065470746100512257520 u^{16}
-4365124819437330950400 u^{15}\]
\vspace{-.4cm}  
\[+272459282567626190070720 u^{14}+19102328814885854400 u^9
+12645306845858008350 u^8\]
\vspace{-.4cm}  
\[ -2606594221714946338575 u^7-48803823903916800 u^2  + 8681150210659989300 \] 
\end{center}
\normalsize
are exactly those numbers of the form $\alpha\beta\gamma$, where 
$(\alpha,\beta,\gamma)$ ranges over all the roots of (\ref{eq:3by3}) in $\C^3$. 
The above 
{\bf univariate reduction} thus tells us that our example indeed has finitely 
many complex roots --- exactly$^{\ref{see}}$ 145, in fact. The above 
polynomial took less than $13$ seconds to compute using a 
naive application of {\bf resultants} and factorization on the 
computer algebra system {\tt Maple}. Interestingly, computing 
the same univariate reduction via a naive application 
of {\bf Gr\"obner bases} (on the same machine with the same version of {\tt 
Maple}) 
takes over $3$ hours and $51$ minutes.$^{\ref{see}}$ 

Admittedly, computing polynomials like the one above can be an 
unwieldy approach to deciding whether (1) has a complex root. An alternative 
algorithm, discovered by Pascal Koiran in \cite{hnam} and improved via theorem 
\ref{main:koi} of section \ref{sec:complex} here, makes a remarkable 
simplification depending on conjectural properties of the distribution of 
prime ideals in number fields.

For instance, an unoptimized implementation of this alternative algorithm 
would run as follows on our example: 
\begin{itemize} 
\item[{\bf Assumption 1}]{ The truth of the {\bf Generalized}\footnote{ The 
{\bf Riemann Hypothesis (RH)} is an 1859 
conjecture equivalent to a sharp quantitative statement on the 
distribution of primes. GRH can be phrased as a generalization of this 
statement to prime ideals in an arbitrary number field. Further background on 
these RH's can be found in \cite{lago,bs}.} {\bf Riemann Hypothesis (GRH)}. } 
\item[{\bf Assumption 2}]{ Access to an {\bf oracle}\footnote{ 
i.e., a machine, or powerful being, which can always instaneously 
and correctly answer such questions. The particular oracle we specify above 
happens to be an {\bf $\np$-oracle} \cite{papa}. } which 
can do the following:  
Given a finite set of polynomials $F\!\subset\!\Z[x,y,z]$ 
and a finite subset $S\!\subset\!\N$, our oracle can decide 
if there is a prime $p\!\in\!S$ such that the mod $p$ reduction of $F$ has a 
root mod in $\Z/p\Z$. } 
\item[{\bf Step 1}]{ Pick a (uniformly distributed) random integer  
$t\!\in\!\{5\cdot 10^6,\ldots,5\cdot 10^6+2\cdot 10^{11}\}$.} 
\item[{\bf Step 2}]{ Using our oracle, decide if there is a prime 
\mbox{$p\!\in\!\{2\cdot 10^{22}\cdot t^3,\ldots,2\cdot 10^{22}\cdot 
(t+1)^3-1\}$} such 
that the mod $p$ reduction of (\ref{eq:3by3}) has a root in $\Z/p\Z$. If so, 
declare  that (\ref{eq:3by3}) has a complex root. Otherwise, declare that 
(\ref{eq:3by3}) has no complex root. \qed } 
\end{itemize} 

The choice of the constants above, 
and the importance of oracle-based algorithms, are detailed 
further in section \ref{sec:complex}. In particular, the constants 
are simply chosen to be large enough to guarantee that, under GRH, the 
algorithm never fails (resp.\ fails with probability $\leq\!\frac{1}{3}$) 
if (\ref{eq:3by3}) has a complex root (resp.\ does not have a 
complex root). Thus, for our example, the algorithm above will 
always give the right answer regardless of the random choice in 
Step 1.  Note also that while the prime we seek above may be quite large, 
the number of {\bf digits} needed to write any such prime is at most 
${\bf 56}$ --- not much bigger than 53, which is the total number 
of digits needed to write down the coefficients and exponent vectors 
of (\ref{eq:3by3}). 
We will explain the complexity-theoretic relevance of this 
fact in section \ref{sec:complex} as well. 
For the sake of completeness, we observe$^{\ref{see}}$  
that the number of real roots of (1) is exactly {\bf 11}. While we will not 
pursue the complexity of real root counting at length in this 
paper, we will quantitatively explore a more general problem over the 
reals. Another example follows. 

\subsection{A Family of Polynomial Inequalities}
\label{sub:real}
In theorem \ref{main:real} of section \ref{sec:real}, we present 
a new bound on the number of connected components of the solution set of any 
collection of polynomial inequalities over the real numbers. Bounds of this 
type have many applications --- for example, 
lower bounds in complexity theory \cite{dl79,sy82} and geometric modelling.

As a simple example, let $S_{a,b}(d,n,p,s)\!\subseteq\!\Rn$ be the solution 
set of the following collection of $p$ equalities and $s$ inequalities: 
\begin{eqnarray}
\label{eq:spike}
\mbox{} \ \ \ \ a_{(\ell,0)}+\left(\sum^{n-1}_{i=1}a_{(\ell,i)}x_i\right)+
\sum^d_{i=1}b_{(\ell,i)}(x_1x_2\cdots x_n)^i & = & 0 \ ; \ \ \ \  
\ell\!\in\!\{1,\ldots,p\} \notag \\
\mbox{} \ \ \ \ a_{(p+\ell,0)}+\left(\sum^{n-1}_{i=1}a_{(p+\ell,i)}x_i\right)+
\sum^d_{i=1}b_{(p+\ell,i)}(x_1x_2\cdots x_n)^i & > & 0 \ ; \ \ \ \ 
\ell\!\in\!\{1,\ldots,s\} 
\end{eqnarray}
for any $d,n,p,s\!\in\!\N$ and real $a_{(i,j)}$ and $b_{(i,j)}$. 

By a bound proved independently by three sets of 
authors between the 1940's and the 1960's \cite{op,milnor,thom}, we 
immediately obtain that $S_{a,b}(d,n,p,s)$ has at most 
$\pmb{(dns+1)(2dns+1)^n}$ connected components.

However, a much sharper bound can be obtained by again looking more closely at 
the monomial term structure involved: 
Let $Q_F$ be the convex hull of the union of the origin $\bO$, the standard 
basis vectors $e_1,\ldots,e_n$ of $\Rn$, and the set of exponent vectors from 
all the polynomials of (\ref{eq:spike}). (In this case, $Q_F$ happens to be 
a bipyramid with one apex at $\bO$ and the other at $(d,\ldots,d)$.) 
Normalizing $n$-dimensional volume, $\vol_n(\cdot)$, so that the 
volume of the $n$-simplex with vertices $\{\bO,e_1,\ldots,e_n\}$ is $1$, 
let $V_F\!:=\vol_n(Q_F)$. Theorem \ref{main:real} then says that
$\min\{n+1,\frac{s+1}{s-1}\}(2s)^nV_F\!=\!\pmb{\min\{n+1,
\frac{s+1}{s-1}\}(2s)^n(d+1)}$ is also an upper bound on the number of 
connected components. 

We have thus improved the older bound by a factor of over $s(dn)^n$ (modulo 
a nonzero multiplicative constant), for this family of 
{\bf semi-algebraic}\footnote{A {\bf semi-algebraic set} is simply a
subset of $\Rn$ defined by the solutions of a finite collection of
polynomial inequalities.} sets. A broader comparison of our bound to earlier 
work appears in section \ref{sub:relatedreal}.  

Let us now fully state our results over $\C$, $\R$, $\Q$, and $\Z$.  

\section{Computing Complex Dimension Faster} 
\label{sec:complex}
Let $f_1,\ldots,f_m\!\in\!\C[x_1,\ldots,x_n]$, $\pmb{F}\!:=\!(f_1,\ldots,
f_m)$, and let $\hn_\C$ denote the problem of deciding whether an 
input $F$ has a complex root.\footnote{We say that $F$ is {\bf feasible} 
(resp.\ {\bf infeasible}) over $\C$ iff $F$ has (resp.\ does not have) a root 
in $\C^n$.}  Also let $\hn$ denote the restriction 
of this problem to polynomials in $\Z[x_1,\ldots,x_n]$. We will 
respectively consider the complexity of $\hn$ and $\hn_\C$ over 
the Turing-machine model and the BSS model over $\C$. 

However, before stating any complexity bounds, let us first clarify our 
notion of input size: With the Turing model, we will assume that any input 
polynomial is given as a sum of monomial terms, with all coefficients {\bf 
and} exponents written in, say, base $2$. The corresponding notion of 
{\bf sparse size} is then simply the total number of bits in all coefficients 
and exponents. For example, the sparse size of $x^D_1+ax^3_1+b$ is $\cO(\log 
D+\log a +\log b)$. The sparse size can be extended to the BSS model 
over $\C$ simply by counting just the total number of bits necessary to write 
down the exponents (thus ignoring the size of the coefficients). 

Note that the number of complex roots of the polynomial 
$x^D_1-1$ is already exponential in its sparse size. 
This behavior is compounded for higher-dimensional polynomial 
systems, and even affects decision problems as well as enumerative problems. 
For example, consider the following theorem. 
\begin{thm}
\cite{plaisted} 
\label{thm:plai} 
$\hn$ is $\np$-hard, even in the special case of two  
polynomial in one variable. More precisely,  
if one can decide whether an arbitrary input polynomial $f\!\in\!\Z[x_1]$
of degree $D$ vanishes at a $D^{\thth}$ root of unity, 
within a number of bit operations polynomial in the sparse size of $f$, 
then $\pp\!=\!\np$. \qed 
\end{thm}
\noindent 
So even for systems such as $f(x_1)\!=\!x^D_1-1\!=\!0$, $\hn$ may be 
impossible to solve within bit complexity polynomial in $\log D$ and the 
sparse size of $f$. An analogue of this result for $\hn_\C$ (theorem 
\ref{thm:smale}) appears in the next section.  

On the other hand, via the classical Sylvester resultant 
\cite[Ch.\ 12]{gkz94} and some basic complexity estimates on arithmetic 
operations \cite{bcs}, it is easy to see that this special case of $\hn$ can 
be decided within a number of bit operations quadratic in $D$ and the sparse 
size of $f$. In complete generality, it is known that $\hn\!\in\!\pspa$ --- 
an important subclass of $\expt$ \cite{koiran}.\footnote{While 
$\pspa$ has important relations to {\bf parallel} algorithms (i.e., 
algorithms where several operations are executed at once by several processors 
\cite{papa}), we will 
concentrate exclusively on {\bf sequential} (i.e., non-parallel) algorithms in 
this paper. } 

Alternatively, if one simply counts arithmetic operations 
(without regard for the size of the intermediate numbers), one can 
similarly obtain an {\bf arithmetic} complexity upper bound of $\cO(D^2)$ for 
the special case of $\hn_\C$ corresponding to the univariate problem mentioned 
in theorem \ref{thm:plai}. 
More generally, it is known that $\hn_\C$ is $\np_\C$-complete\footnote{This 
is the analogue of $\np$-complete for the BSS model over $\C$ \cite{bcss}.} 
\cite{bss,shub}. 

Curiously, efficient {\bf randomization-free} algorithms for  
$\hn$ and $\hn_\C$ are hard to find in the literature. So we present such an 
algorithm, with an explicit complexity bound, for a problem including 
$\hn_\C$ as a special case.  
\begin{main}
\label{main:complex}
Let $Z_F$ be the zero set of $F$ in $\Cn$ and $\dim Z_F$ the complex dimension 
of $Z_F$. Also let $\bO$ be the origin, and $e_1,\ldots,e_n$ the standard 
basis vectors, in $\Rn$. Normalize $n$-dimensional volume $\vol_n(\cdot)$ so 
that the volume of the standard $n$-simplex (with vertices 
$\bO,e_1,\ldots,e_n$) is $1$. Finally, let $k$ be the total number of monomial 
terms in $F$ (counting repetitions between distinct $f_i$) and let 
$Q_F$ be the convex hull of the union of $\{\bO,e_1,\ldots,e_n\}$ and the set 
of all exponent vectors of $F$. Then there is a 
deterministic\footnote{i.e., randomization-free} 
algorithm which computes $\dim Z_F$, and thus solves $\hn_\C$, within 
$\cO(n^4kM^{2.376}_FV^5_F+nk\log(m+n))$ arithmetic operations, where 
$V_F\!:=\!\vol_n(Q_F)$ and $M_F$ is no larger than the maximum number of 
lattice points in any translate of $(n+1)Q_F$. 
\end{main}
\noindent 
Via a height\footnote{The (absolute multiplicative) {\bf height} of an 
algebraic number $\zeta$ is an important number-theoretic invariant 
related to the minimal polynomial of $\zeta$ over $\Z$. Height 
bounds are also intimately related to more pedestrian quantities 
like the maximum absolute value of a coordinate of an isolated root of a 
polynomial system, so we use the term ``height'' in this collective sense. 
Further details on heights, 
and their extension to $\Cn$, can be found in \cite{sil,gregoheight,cool}. } 
estimate from theorem \ref{main:height} later in 
this section one can also derive a similar bound on the bit complexity of 
$\hn$.  We clarify the benefits of our result over earlier bounds  
in section \ref{sub:relatedcomplex}. The algorithm for theorem 
\ref{main:complex}, and its correctness proof, are 
stated in section \ref{sub:proofcomplex}. The techniques 
involved will also be revisited in our discussion of quantifier 
prefixes over $\Z$ in section \ref{sec:int}.  

There is, however, a fundamentally different approach which, given the truth 
of GRH, places $\hn$ in an even better complexity class. 
First recall that randomized decision algorithms which answer incorrectly 
with probability, say, $\leq\!\frac{1}{3}$, and for which the number of bit 
operations 
and random bits needed is always polynomial in the input size, define the 
complexity class $\bpp$.\footnote{We emphasize that such algorithms can give 
different answers when run many times on the same input. However, by accepting 
the most popular answer of a large sample, the error probability can be made 
arbitrarily small.} Recall also that when a $\bpp$ algorithm is augmented by 
an oracle in $\np$, and the number of 
oracle-destined bits is always polynomial in the input size, one obtains the 
class $\bpp^\np$. Finally, when just {\bf one} oracle call is allowed 
in a $\bpp^\np$ algorithm, one obtains the {\bf Arthur-Merlin 
class} $\am$ \cite{zachos}. 
\begin{thm}
\label{thm:koi}
\cite{hnam}
Assuming the truth of GRH, $\hn\!\in\!\am$. \qed 
\end{thm} 
\noindent 
While probabilistic algorithms for $\hn$ (and more general problems) 
have certainly existed at least since the early 1980's, the above theorem 
is the first and only example of an algorithm for $\hn$ requiring a number 
of bit operations just {\bf polynomial} in the input size, albeit 
modulo two strong assumptions. 

In view of the vast literature on GRH from both number theory 
and theoretical computer science, the study of algorithms depending on GRH is 
not unreasonable. For example, the truth of GRH implies a polynomial-time 
algorithm for deciding whether an input integer is prime \cite{miller}. 
Likewise, in view of the continuing open status of the 
$\pp\stackrel{?}{=}\np$ question, oracle-based results are well-accepted 
within theoretical computer science.\footnote{It turns out that $\pp\!=\!\np$ 
also implies the existence of a polynomial-time algorithm for primality 
testing \cite{pratt}.} 
In particular, 
Koiran's conditional result gives the smallest complexity class known to 
contain $\hn$. Indeed, independent of GRH, while 
it is known that $\np\!\subseteq\!\am\!\subseteq\!\pspa$ \cite{papa}, the 
properness of each inclusion is still an open problem.

The simplest summary of Koiran's algorithm is that it uses 
reduction modulo specially selected primes to decide feasibility over $\C$. 
(His algorithm is unique in this respect since all previous 
algorithms for $\hn$ worked primarily in the ring $\C[x_1,\ldots,x_n]/\langle 
F \rangle$.) The key observation behind Koiran's algorithm is that 
an $F$ infeasible (resp.\ feasible) over $\C$ will have roots in $\Z/p\Z$ 
for only finitely many (resp.\ a positive density of) primes $p$. 

A refined characterization of the difference between positive and zero 
density can be given in terms of our framework as follows:  
\begin{thm} 
\label{main:koi}
Following the notation above, assume now that 
$f_1,\ldots,f_m\!\in\!\Z[x_1,\ldots,x_n]$, let\footnote{We point out 
that in \cite{hnam}, the notation $\sigma(F)$ was instead used for 
a different quantity akin to $2+mD$. } $\sigma(F)$ be the maximum 
of $\log|c|$ as $c$ ranges over the coefficients of all the monomial terms of 
$F$, and let $D$ be the maximum total degree of any $f_i$. Then there exist 
$a_F,A_F\!\in\!\N$, with the following properties: 
\begin{itemize}
\item[(a)]{$F$ infeasible over $\C \Longrightarrow$ the reduction of 
$F$ mod $p$ has a root in $\Z/p\Z$ for at most $a_F$ distinct primes $p$, 
and $a_F\!=\!\cO(n^3DV_F(4^nD\log D + \sigma(F)+\log m))$.}  
\item[(b)]{Given the truth of GRH, $F$ feasible over $\C \Longrightarrow$ for 
each $t\!\geq\!4963041$, the sequence 
\mbox{$\{A_Ft^3,\ldots,A_F(t+1)^3-1\}$} contains 
a prime $p$ such that the reduction of $F$ mod $p$ has a root in $\Z/p\Z$. 
Furthermore, we can take 
$A_F\!=\!O\left([V_F\sigma(h_F)(n\log D+\log\sigma(F))]^2\right)$, where $h_F$ 
is the polynomial defined in theorem \ref{main:height} below. } 
\end{itemize}
In particular, the bit-sizes of $a_F$ and $A_F$ are both 
\mbox{$\cO(n\log D+\log \sigma(F))$} --- 
sub-quadratic in the sparse size of $F$. Simple explicit 
formulae for $a_F$ and $A_F$ appear in remarks \ref{rem:shebanga} and 
\ref{rem:shebangb} of section \ref{sub:proofcomplex}. 
\end{thm}

Via theorem \ref{main:koi}, Koiran's algorithm for $\hn$ can be 
paraphrased as follows:\footnote{We point out that, to the best of the 
author's knowledge, this is the first time that the constants 
underlying Koiran's algorithm have been made explicit.}  
\begin{itemize} 
\item[{\bf Assumption 1}]{ The truth of GRH.} 
\item[{\bf Assumption 2}]{ Access to an $\np$-oracle. } 
\item[{\bf Step 1}]{ Pick a (uniformly distributed) random integer  
\mbox{$t\!\in\!\{4963041,\ldots,4963041+3a_F\}$.} } 
\item[{\bf Step 2}]{ Using our oracle, decide if there is a prime 
\mbox{$p\!\in\!\{A_Ft^3,\ldots,A_F(t+1)^3-1  
\}$} such that $F$ has a root mod $p$. If so, declare  
that $F$ has a complex root. Otherwise, declare that $F$ has no complex root. 
\qed } 
\end{itemize} 
\noindent 
In particular, it follows immediately from theorem \ref{main:koi} that 
the algorithm above is indeed an $\am$ algorithm, and that the error 
probability is $\leq\!\frac{1}{3}$. Better still, the error probability can be 
replaced by an arbitrarily 
small constant $\eps$ (keeping the same asymptotic complexity), simply  
by replacing $3a_F$ by $\frac{1}{\eps}a_F$ in Step 1 above. 

The proof of theorem \ref{main:koi} is based in part on a particularly 
effective form of univariate reduction. 
\begin{thm} 
\label{main:height}
Following the notation above, and the assumptions of theorem \ref{main:koi}, 
there exist a univariate polynomial $h_F\!\in\!\Z[u_0]$ and 
a point $u_F\!:=\!(u_1,\ldots,u_n)\!\in\Zn$ with the following 
properties: 
\begin{enumerate} 
\setcounter{enumi}{-1}  
\item{The degree of $h_F$ is $\leq\!V_F$.} 
\item{For any irreducible component $W$ of $Z_F$, there is a 
point $(\zeta_1,\ldots,\zeta_n)\!\in\!W$ such that 
$u_1\zeta_1+\cdots+u_n\zeta_n$ is a root of $h_F$. Conversely, 
if $m\!\leq\!n$, all roots of $h_F$ arise this way. } 
\item{$F$ has only finitely many complex roots $\Longrightarrow$ the splitting 
field of $h_F$ over $\Q$ is exactly the field $\Q[x_i \; \; | \; \; 
(x_1,\ldots,x_n)\!\in\!\Cn \text{ \ is \ a \ root \ of \ } F]$. } 
\item{The coefficients of $h_F$ satisfy 
$\sigma(h_F)\!=\!\cO\left(M_F[\sigma(F)+m(n\log D+\log m)]+n^2V_F\log 
D\right)$ and, when $m\!\leq\!n$, 
$\sigma(h_F)\!=\!\cO(M_F\sigma(F)+n^2V_F\log D)$.}  
\item{$m\!\leq\!n \Longrightarrow$ the deterministic arithmetic complexity 
of computing $u_F$, and all the coefficients of $h_F$, is 
$\cO(n^3M^{2.376}_F V^5_F)$. } 
\item{We have $\log(1+|u_i|)\!=\!\cO(n^2\log D)$ for all $i$.} 
\end{enumerate} 
\end{thm} 
\noindent 
Note that we have thus obtained the existence of points of 
bounded height on the positive-dimensional part of $Z_F$, 
as well as a bound on the height of any point in the 
zero-dimensional part of $Z_F$. Put more simply, via a slight 
variation of the proof of theorem \ref{main:height}, we obtain 
the following useful bound:  
\begin{thm} 
\label{main:size} 
Following the notation of theorem \ref{main:height}, 
any irreducible component $W$ of $Z_F$ contains a point 
$(x_1,\ldots,x_n)$ such that for all $i$, either 
$x_i\!=\!0$ or\\ 
$|\log|x_i||\!=\!\cO\left(M_F[\sigma(F)+m(n\log D+\log m)]\right)$. 
Furthermore, when $m\!\leq\!n$, the last upper bound can be improved 
to $\cO(M_F\sigma(F))$. 
\qed 
\end{thm}

Our final result over $\C$ is a refinement of theorem \ref{main:height} which 
will help simplify the proofs of our results in section \ref{sec:int} on 
integral points. 
\begin{thm} 
\cite{gcp}
\label{main:unired} 
Following the notation of theorem \ref{main:height}, 
one can pick $u_F$ and $h_F$ (still satisfying (0)--(5)) so that there 
exist $a_1,\ldots,a_n\!\in\!\N$ and $h_1,\ldots,h_n\!\in\!\Z[u_0]$ with the 
following properties: 
\begin{enumerate} 
\setcounter{enumi}{5}
\item{ The degrees of $h_1,\ldots,h_n$ are all bounded above 
by $V_F$.} 
\item{ For any root $\theta\!=\!u_1\zeta_1+\cdots+u_n\zeta_n$ of $h_F$, 
$\frac{h_i(\theta)}{a_i}\!=\!\zeta_i$ for all $i$. } 
\item{For all $i$, both $\log a_i$ and 
$\sigma(h_i)$ are bounded above by $\cO(V^2_F\sigma(h_F))$.}  
\item{$m\!\leq\!n \Longrightarrow$ the deterministic arithmetic complexity of 
computing all the coefficients of $h_1,\ldots,h_n$ is 
$\cO(n^4M^{2.376}_F V^5_F)$. } 
\end{enumerate}
\end{thm} 

Explicit formulae for all these asymptotic estimates, as well as 
their proofs, appear in remarks \ref{rem:height}, \ref{rem:size}, and 
\ref{rem:denom} of section \ref{sub:proofcomplex}. However, let us first 
compare these quantitative results to earlier work. 

\subsection{Related Results Over $\C$} 
\label{sub:relatedcomplex}  
Solving $\hn_\C$ too quickly also leads to unexpected collapses of 
complexity classes as follows. 
\begin{thm}
\label{thm:smale} 
Suppose there is an algorithm (on a BSS machine over $\C$) which decides  
whether an arbitrary input polynomial $f\!\in\!\C[x_1]$ of degree $D$ 
vanishes at a $D^{\thth}$ root of unity, within a number of arithmetic 
operations polynomial in the sparse size of $f$. Then $\np\!\subseteq\!\bpp$. 
\qed 
\end{thm}
\noindent 
This result is originally due to Steve Smale and a proof appears 
in \cite{real}. It is currently believed that the inclusion 
$\np\!\subseteq\!\bpp$ is quite unlikely. 

Curiously, finding (as opposed to deciding the existence of) roots 
for even a seemingly innocent univariate polynomial can lead to 
undecidability in the BSS model over $\C$: 
\begin{thm} 
Determining whether an {\bf arbitrary} $x_0\!\in\!\C$ converges to a 
root of $x^3-2x+2\!=\!0$ under Newton's method is undecidable, relative to 
the BSS model over $\C$. \qed 
\end{thm} 
\noindent 
This result follows easily via a dynamics result of Barna \cite{barna} 
and the proof appears in \cite[Sec.\ 2.4]{bcss}. One should of 
course note that this result in no way prevents one from finding 
{\bf some} $x_0$ which converges to a root of $x^3-2x+2$. So this 
result is a more a reflection of the subtlety of dynamics than 
the limits of the BSS model. 

As for the other results of section \ref{sec:complex}, we point out that we 
have tried to balance generality, sharpness, and ease of proof in our 
bounds. In particular, our bounds fill a lacuna in the literature where 
earlier bounds seemed to sacrifice generality for sharpness, or vice-versa. 

To clarify this trade-off, first note that 
$\cI_F\!\leq\!V_F\!\leq\!D^n$, 
where $\cI_F$ is the number of irreducible components of $Z_F$. 
(The first inequality follows immediately from theorem 
\ref{main:height}, while the 
second follows from the observation that $Q_F$ always lies in a 
copy of the standard $n$-simplex scaled by a factor of $D$.) 
So depending on the shape of $Q_F$, and thus 
somewhat on the sparsity of $F$, one can typically expect $V_F$ 
to be much smaller than $D^n$. For example, our $3\times 3$ system 
from section \ref{sub:3by3} gives $D^n\!=\!13824$ and $V_F\!=\!243$. 
Setting $p\!=\!n$ and $s\!=\!0$ in the example from section \ref{sub:real}, 
it is easy to see that the factor of improvement can even reach 
$D^{n-1}$, if not more. 

As for the quantities $k$ and $M_F$, we will see 
in lemma \ref{lemma:respert} of section \ref{sub:res} that 
$k\!\leq\!m(V_F+n)$ and $M_F\!\leq\!\begin{pmatrix} 
nD+1\\ n\end{pmatrix}\!=\!\cO(e^n(nD+1)^n)$. Furthermore, just as 
$V_F$ is a much more desirable complexity measure than $D^n$, we point out 
that the preceding bound on $M_F$ is frequently overly pessimistic: for 
example, 
$M_F\!=\!\cO(V_F)$ for fixed $n$. The true 
definition of $M_F$ appears in section \ref{sub:res}.  

Our algorithm for computing $\dim Z_F$ thus gives the first deterministic 
complexity bound which is polynomial in $V_F$ and $M_F$. In particular, 
while harder problems were already known to admit $\pspa$ complexity bounds, 
the corresponding complexity bounds were either polynomial (or worse) in 
$D^n$, or stated in terms of a non-uniform computational model.\footnote{ 
For example, some algorithms in the literature are stated in terms 
of {\bf arithmetic networks}, where the construction of the 
underlying network is not included in the complexity estimate. }  
Our algorithm for the computation of $\dim Z_F$ thus gives a significant 
speed-up over earlier work. 

For example, via the work of Chistov and Grigoriev from the early 1980's on 
quantifier elimination over $\C$ \cite{chigo}, it is not hard to derive a 
deterministic arithmetic complexity 
bound of $\cO((mD)^{n^4})$ for the computation of $\dim Z$. More recently, 
\cite{giustiheintz} gave a randomized arithmetic complexity bound of 
$m^{\cO(1)}D^{\cO(n)}$. Theorem \ref{main:complex} thus clearly 
improves the former bound. Comparison with the latter bound is 
a bit more difficult since the exponential constants and derandomization 
complexity are not explicit in \cite{giustiheintz}. 

As for faster algorithms, one can seek complexity bounds which are polynomial 
in even smaller quantities. For example, if one has an irreducible algebraic 
variety $V\!\subseteq\!\Cn$ of complex dimension $d$, one can 
define its {\bf affine geometric degree}, $\delta(V)$, to be the number of 
points in $V\cap H$ where $H$ is a generic $(n-d)$-flat.\footnote{ We explain 
the term ``generic'' in sections \ref{sec:int} and \ref{sub:h3}.} More 
generally, we can 
define $\delta(Z_F)$ to be the sum of $\delta(V)$ as $V$ ranges over all 
irreducible components of $Z_F$. It then follows (from theorem 
\ref{main:complex} and a consideration of intersection multiplicities) 
that $\cI_F\!\leq\!\delta(Z_F)\!\leq\!V_F$. 
Similarly, one can attempt to use mixed volumes of several polytopes (instead 
of a single polytope volume) to lower our bounds. 

We have avoided refinements of this nature for the sake of simplicity. 
Another reason it is convenient to have bounds in terms of $V_F$ 
is that the computation of $\delta(Z_F)$ is even more subtle 
than the computation of polytopal $n$-volume. For example, when $n$ is 
fixed, $\vol_n(Q)$ can be computed in polynomial time simply by 
triangulating the polytope $Q$ and adding together the volumes of the 
resulting $n$-simplices 
\cite{volcomplex}. However, merely deciding $\delta(Z_F)\!>\!0$ is 
already $\np$-hard for $(m,n)\!=\!(2,1)$, via theorem \ref{thm:plai}. 
As for varying $n$, computing 
$\delta(Z_F)$ is $\#\pp$-hard, while the computation of polytope volumes is 
$\#\pp$-complete.\footnote{ $\#\pp$ is the analogue of $\np$ for enumerative 
problems (as opposed to decision problems) \cite{papa}. } 
(The latter result is covered in \cite{volcomplex,kls}, while the former 
result follows immediately from the fact that the computation of $\delta(Z_F)$ 
includes the computation of $V_F$ as a special case.) More practically, for 
any fixed $\eps_1,\eps_2\!>\!0$, there 
is an algorithm which runs in time polynomial in the sparse encoding of $F$ 
(and thus polynomial in $n$) which produces a random variable that is within a 
factor of $1-\eps_1$ of $\vol_n(Q_F)$ with probability $1-\eps_2$ \cite{kls}. 
The analogous result for mixed volume is known only for certain families 
of polytopes \cite{gs00}, and the existence of such a result for 
$\delta(Z_F)$ is still an open problem. 

In any event, we point out that improvements in terms of $\delta(Z_F)$ 
for our bounds are possible, and these will be pursued in a forthcoming  
paper. Similarly, the exponents in our complexity bounds can 
be considerably lowered if randomization is allowed. Furthermore, Lecerf has 
recently announced a randomized arithmetic complexity 
bound for computing $\dim Z_F$ which is polynomial in 
$\max_i\{\delta(Z_{(f_1,\ldots,f_i)})\}$ 
\cite{lecerf}.\footnote{ The paper \cite{lecerf} actually solves the  
harder problem of computing an algebraic description of a non-empty 
set of points in every irreducible component of $Z_F$, and distinguishing 
which component each set belongs to.} However, the complexity of 
derandomizing Lecerf's algorithm is not yet clear. 

As for our result on prime densities (theorem \ref{main:koi}), part (a) 
presents the best current bound polynomial in $V_F$ and $M_F$. An earlier 
density bound, polynomial in $D^{n^{\cO(1)}}$ instead, 
appeared in \cite{hnam}. 

Part (b) of theorem \ref{main:koi} appears to be new, and 
makes explicit an allusion of Koiran in \cite{hnam}. 
\begin{rem}
\label{rem:foo} 
We point out that we cheated slightly in our refinement of 
Koiran's algorithm: We did not take the complexity of computing 
$V_F$ into account. (It is easy to see that this is what dominates 
the randomized bit complexity of the algorithm.) This can be corrected, and 
perhaps the simplest way is to 
replace every occurence of $V_F$ with $D^n$ in our bounds for $M_F$, $a_F$, 
and $A_F$. Alternatively, if one want to preserve polynomiality in $V_F$, 
one can instead apply the polynomial-time randomized approximation techniques 
of \cite{kls} to $V_F$, and make a minor adjustment to the error 
probabilities. \qed 
\end{rem} 
\begin{rem} 
Pascal Koiran has also given an $\am$ algorithm 
(again depending on GRH) for deciding whether the complex dimension 
of an algebraic set is less than some input constant \cite{koiran}. 
\qed 
\end{rem} 

Regarding our height bound, the only other results stated in polytopal terms 
are an earlier version of theorem \ref{main:height} announced in 
\cite{stoc99}, and independently 
discovered bounds in \cite[Prop.\ 2.11]{cool} and \cite[Cor.\ 
8.2.3]{maillot}. The bound from \cite{cool} applies to a slightly 
different problem, but implies (by intersecting with a generic 
linear subspace with reasonably bounded coefficients)\footnote{ 
Martin Sombra pointed this out in an e-mail to the author.} a bound of  
$\cO((4^nD\log n + n\sigma(F))V_F)$ for our setting. Furthermore, 
by examining a key ingredient in their proof (Proposition 1.7 from 
\cite{cool}), their bound can actually be improved to 
$\cO(DM_F\log n+nV_F\sigma(F))$. The last bound is thus close to ours, 
and can be better when $m$ and $\sigma(F)$ are large and $n$ is small. 
The bound from \cite[Cor.\ 8.2.3]{maillot} uses 
Arakelov intersection theory, holds only for $m\!=\!n$, and the statement is 
more intricate (involving a sum of several mixed volumes). So it is not yet 
clear when \cite[Cor.\ 8.2.3]{maillot} is better than theorem 
\ref{main:height}. In any case, our result has a considerably simpler 
proof than either of these two alternative bounds: We use only resultants and 
elementary linear algebra and factoring estimates.  

We also point out that the only earlier 
bounds which may be competitive with theorems \ref{main:height} and 
\ref{main:size}, \cite[Prop.\ 2.11]{cool}, and \cite[Cor.\ 8.2.3]{maillot} 
are polynomial in 
$e^n(nD+1)^n$ and make various non-degeneracy hypothesis, e.g., 
$m\!=\!n$ and no singularities for $Z_F$ (see \cite{cannyphd} and 
\cite[Thm.\ 5]{gregogap}). As for 
bounds with greater generality, the results of \cite{fgm} imply a height bound 
for general quantifier elimination which, unfortunately, has a factor of the 
form $2^{(n\log D)^{\cO(r)}}$ where $r$ is the number of 
quantifier alternations \cite{hnam}. 

As for theorem \ref{main:unired}, 
the approach of rational univariate representations ({\bf RUR}) for the 
roots of polynomial systems dates back to Kronecker. RUR also goes under the 
name of ``effective primitive element theorem'' and important precursors to 
theorem \ref{main:unired}, with respective complexity bounds polynomial 
in $e^n(nD+1)^n$ and $D^{n^{\cO(1)}}$, are stated in \cite{pspace} and 
\cite[Thm.\ 4]{hnam}. Nevertheless, the use of 
{\bf toric resultants} (cf.\ section \ref{sub:proofcomplex}), 
which form the core of our algorithms here, was 
not studied in the context of RUR until the late 1990's (see, e.g., 
\cite{gcp}). In particular, theorem \ref{main:unired} appears to be the 
first statement giving bounds on $\sigma(h_i)$ which are polynomial 
in $V_F$. As for computing $h,h_1,\ldots,h_n$ faster, an algorithm for 
RUR with randomized complexity polynomial in 
$\max_i\{\delta(Z_{(f_1,\ldots,f_i)})\}$ was derived in 
\cite{gls99}. However, their 
algorithm makes various nondegeneracy assumptions (such as $m\!=\!n$ and that 
$F$ form a complete intersection) and the derandomization 
complexity is not stated. 

The remaining bottle-neck in improving our complexity and 
height bounds stems from the exponentiality in $n$ present in the quantity 
$M_F$. However, the resulting exponential factor,  
which is currently known to be at worst $\cO(e^n)$ (cf.\ 
lemma \ref{lemma:respert} of section \ref{sub:res}), can be reduced to 
$\cO(n)$ in certain cases. 
In general, this can be done whenever there exists an expression for a 
particular toric resultant (cf.\ section \ref{sub:proofcomplex}) as a single 
determinant, or the divisor of a determinant, of a matrix of size $\cO(nV_F)$. 
The existence of such formulae has been proven in various cases, e.g., when 
all the Newton polytopes are axis-parallel parallelepipeds \cite{wz}. Also, 
such formulae have been observed (and constructed) experimentally in various 
additional cases of practical interest \cite{emican}. Finding compact formulae 
for resultants is an area of active research which thus has deep implications 
for the complexity of algebraic geometry.  

Finally, we note that 
we have avoided Gr\"obner basis techniques because there are currently 
no known complexity or height bounds polynomial in $V_F$ (or even 
$M_F$) using Gr\"obner bases for the problems we consider. A further 
complication is that 
there are examples of ideals, generated by polynomials of degree $\leq\!5$ in 
$\cO(n)$ variables, where every Gr\"obner basis has a generator of degree 
$2^{2^n}$ \cite{mm}. This is one obstruction to deriving sharp explicit 
complexity bounds via a naive application of Gr\"obner bases. Nevertheless, 
we point out that Gr\"obner bases are well-suited for other difficult 
algebraic problems, and their complexity is also an area of active research. 
  
\section{Polytope Volumes and Counting Pieces of Semi-Algebraic Sets} 
\label{sec:real}
Continuing our theme of measuring algebraic-geometric complexity in 
combinatorial terms, we will see how to bound the number of connected 
components of a semi-algebraic set in terms of polytope volumes. However, let 
us first see an unusual example of how input encoding influences 
computational complexity, as well as geometric complexity, over the real 
numbers.  

Recall that a {\bf straight-line program (SLP)} presents 
a polynomial as a sequence of subtractions and multiplications, starting from 
a small set of constants and variables \cite{bcs,bcss}. (Usually, the only 
constant given a priori is $1$.) The {\bf SLP size} of a polynomial 
$f\!\in\!\Z[x_1,\ldots,x_n]$ is then just the minimum of the total 
number of operations needed by any SLP evaluating to $f$. Thus, while 
$(x+2^{2^2})^{1000}-2^{2^{2^3}}$ has a large sparse size, its SLP size is 
easily seen to be quite small, via standard recursive tricks such as repeated 
squaring.  SLP's are thus a more powerful encoding than the sparse encoding, 
since the SLP size of a polynomial is trivially bounded from above by a linear 
function of its sparse size.  

Consider the following corollary of theorem \ref{thm:plai}. 
\begin{cor} 
\label{cor:realuni} 
If one can decide whether an arbitrary $f\!\in\!\Z[x_1]$ has a real root, 
within a number of bit operations polynomial in the SLP size of $f$, 
then $\pp\!=\!\np$. \qed 
\end{cor}
\noindent
Thus the hardness of feasibility testing we've observed earlier over $\C$ 
persists over $\R$, albeit relative to a smaller complexity measure. 
Peter B\"{u}rgisser observed the following simple proof of 
this corollary in 1998: Assuming the hypothesis above, consider 
the polynomial system $G\!:=\!(f(w),w(z+i)-iz)$. 
Then $f$ has a real root $\Longleftrightarrow 
G$ has a root $(w,z)$ with $w$ on the unit circle, and our assumption 
thus implies the existence of a polynomial-time algorithm (relative now to the 
SLP encoding) for detecting whether certain systems of two polynomials 
in two variables have a root $(w,z)$ with $w$ on the unit circle. 
This in turn implies an algorithm, requiring a number 
of bit operations just polynomial in the sparse size of $f$, for deciding if 
a univariate polynomial $f$ has a root on the unit circle. This is not quite 
the same problem as the special case of $\hn$ from theorem \ref{thm:plai}, 
but it is nevertheless known to be $\np$-hard as well \cite{plaisted}. 
So we finally obtain $\pp\!=\!\np$ from our initial assumption and our 
corollary is thus proved. 

Another complication with detecting the existence of real roots too quickly 
is that the number of real roots, even for a single univariate polynomial, 
can be exponential in the SLP size. (This fact is {\bf not} implied by 
our earlier example of $x^D_1-1$.) To see why, simply consider the recursion  
$g_{j+1}\!:=\!4g_j(1-g_j)$ with $g_1\!:=\!4x(1-x)$. It is then easily 
checked\footnote{ This example is well-known in dynamical systems, and the 
author thanks Gregorio Malajovich for pointing it out. } 
that $g_j(x)-x$ has $2^j$ roots in the open interval $(0,1)$, but an SLP size 
of just $\cO(j)$.  

It is an open question whether corollary \ref{cor:realuni}  
holds relative to {\bf sparse} size. More to the point, 
the influence of sparse size on the number of 
{\bf real} roots of polynomial systems remains a deep 
open question. For instance, the classical {\bf Descartes rule of signs} 
states that any univariate polynomial with real coefficients 
and $k$ monomial terms has at most $2k+1$ real roots. However, 
the best known bounds on the number of isolated real roots 
for $2$ polynomials in $2$ unknowns are already exponential 
in the number of monomial terms, even if one restricts to roots 
with all coordinates positive (cf.\ section \ref{sub:relatedreal}).  

However, one can at least give bounds which are linear  
in a suitable polytope volume, which apply even in the  
the more general context of polynomial inequalities. 
\begin{main}
\cite{real}
\label{main:real} 
Let $f_1,\ldots,f_{p+s}\!\in\!\R[x_1,\ldots,x_n]$ and suppose 
$S\!\subseteq\Rn$ is the solution set of the following 
collection of polynomial inequalities: 
\begin{eqnarray*} 
f_i(x)\!&\!=\!&\!0, \ \ \ i\!\in\!\{1,\ldots,p\} \\ 
f_{p+i}(x)\!&\!>\!&\!0, \ \ \ i\!\in\!\{1,\ldots,s\} 
\end{eqnarray*} 
Let $Q_F\!\subset\!\Rn$ be the convex hull of the 
union of $\{\bO,\hat{e}_1,\ldots,\hat{e}_n\}$ 
and the set of all $a$ with $x^a\!:=\!x^{a_1}_1\cdots x^{a_n}_n$ a  
monomial term of some $f_i$. Then $S$ has at most 
\[ \mbox{}\hspace{-.3cm}\min\{n+1,\frac{s+1}{s-1}\}2^ns^nV_F \ 
(\mathrm{for \ } s\!>\!0) \ \ \mathrm{or}   \ \ 2^{n-1}V_F \ 
(\mathrm{for \ } s\!=\!0) \]  
connected components, where $V_F\!:=\!\vol_n(Q_F)$. \qed 
\end{main} 

In closing this brief excursion into semi-algebraic geometry, we point out 
that unlike the complex case, it is not yet known  
whether $V_F$ is an upper bound on the number of {\bf real} connected 
components. This is because a complex component may contribute two or more 
real connected components. Nevertheless, it is quite possible that the factors 
exponential in $n$ in our bounds may be removed from our bounds in the near 
future. 

\subsection{Related Results Over $\R$} 
\label{sub:relatedreal}
We first recall the following important result 
relating sparse size and real roots for certain non-degenerate 
polynomial systems. (Recall also that the {\bf positive orthant} 
of $\Rn$ is the subset $\{(x_1,\ldots,x_n) \; | \; x_i\!>\!0 
{\text \ for \ all \ } i\}$.)  
\begin{kho} 
{\bf (Special Case)}\footnote{Khovanski's Theorem on Fewnomials 
actually holds for a more general class of functions --- the so-called 
{\bf Pfaffian} functions \cite{few}. } \cite[Sec.\ 3.12, Cor.\ 6]{few}
Following the notation of theorem \ref{main:real}, suppose $p\!=\!n$,  
$s\!=\!0$, and the Jacobian matrix of $F$ is invertible at any  
complex root of $F$. Also let $k'$ be the number of exponent 
vectors which appear in at least one of $f_1,\ldots,f_n$. Then $F$ has at 
most $(n+1)^{k'} 2^{k'(k'-1)/2}$ real roots in the positive orthant. \qed 
\end{kho} 
\noindent 
For example, Khovanski's bound readily implies that our $3\times 3$ example 
from section \ref{sub:3by3} has at most 
$8\cdot 4^9\cdot 2^{36}\!=\!\mathbf{144115188075855872}$ real roots --- quite 
a bit more than $972$ (the estimate from theorem \ref{main:real} above) or 
$11$ (the true number of real roots). Nevertheless, we emphasize that his 
theorem was a major advance, giving the first bound on the number of real 
roots independent of the degree of the input polynomials. 

As for other more general results, Khovanski also gave  
bounds on the {\bf Betti numbers}\footnote{These are more subtle 
cohomological invariants which include the number of connected components 
as a special case (see, e.g., \cite{munkres} for further details).} of 
non-degenerate real algebraic varieties \cite[Sec.\ 3.14, Cor.\ 5]{few}. 
Similarly, these results (which thus require $p\!\leq\!n$ and $s\!=\!0$) 
become more practical as the polynomial degrees 
grows and the number of monomial terms remains small. 

Closer to our approach, Benedetti, Loeser, and Risler 
independently derived a polytopal upper bound on the number connected 
components of a real algebraic variety in \cite[Prop.\ 3.6]{blr}. 
Their result, while applying only in the case where $p\!\leq\!n$ and $s\!=\!0$, 
can give a better bound when the number of equations $p$ is a small constant 
and $n$ is large. We also point out that their result has a more complicated 
statement than ours, involving a recursion in terms of mixed volumes of 
projections of polytopes. 

The only other known bounds on the number of 
connected components appear to be linear in $D^n$. For example, 
a bound derived by Oleinik, Petrovsky, Milnor, and Thom 
before the mid-1960's \cite{op,milnor,thom} gives $D(2D-1)^{n-1}$ for 
$s\!=\!0$ and $(sD+1)(2sD+1)^n$ for $s\!>\!0$. An improvement, also 
polynomial in $D^n$, was given recently by Basu \cite{basu}:  
\mbox{$(p+s)^n\cO(D)^n$,} where the implied constant is not stated explicitly.  
For $s\!>\!0$ our bound is no worse than 
$\min\{n+1,\frac{s+1}{s-1}\}(2sD)^n$ --- better than both preceding bounds 
and frequently much better. For $s\!=\!0$ our bound is no worse 
than $2^{n-1}D^n$ --- negligibly worse than the oldest bound, but 
asymptotically better than Basu's bound.  

For the sake of brevity, we have mainly focused on one combinatorial aspect of 
semi-algebraic sets. So let us at least mention a few additional 
complexity-theoretic 
references: Foundational results on the complexity of solving (or counting 
the roots of) polynomial systems over $\R$ can be 
found in \cite{marie}, and faster recent algorithms can be found in 
\cite{esa,moupan}.  More generally, there are 
algorithms known for quantifier elimination over any real closed field 
\cite{renegar,cannyquant,bpr}. 

Curiously, the best current complexity bounds for the problems over $\R$ just 
mentioned are essentially the same as those for the corresponding problems 
over $\C$. Notable recent exceptions include \cite{bank} and \cite{rojasye} 
where the complexity bounds depending mainly on quantities relating only to 
the underlying real geometry. (The first paper deals with finding a 
point in every connected component of a semi-algebraic set, while the second 
paper deals with approximating the real roots of a trinomial within time 
quadratic in $\log D$.) Also, with the exception of 
\cite{bank,esa,moupan,rojasye}, all the preceding references present 
complexity bounds depending on $n$ and $D^n$, with no 
mention of sharper quantities like $V_F$. 

An interesting question which remains is whether feasibility over $\R$ can 
be decided within the {\bf polynomial hierarchy} (a collection of 
complexity classes suspected to lie below $\pspa$ \cite{papa}), with or 
without GRH. As we will see now, this can be done 
over $\Q$ (at least in a restricted sense) as well as $\C$. 

\section{The Generalized Riemann Hypothesis and Detecting 
Rational Points} 
\label{sec:rat} 
Here we will return to considering computational complexity 
estimates: We show that deciding feasibility over $\Q$, for 
most polynomial systems, lies within the polynomial hierarchy, assuming 
GRH. To fix ideas, let us begin with the case of a single univariate 
polynomial. 
\begin{thm}
\cite{lenstra}
\label{thm:lenstra} 
Suppose $f\!\in\!\Z[x_1]$ and $\pm\frac{p}{q}\!\in\!\Q$ is a root of $f$, 
with $p,q\!\in\!\N$ and $\gcd(p,q)\!=\!1$. Then $\log p$, $\log q$, 
and the number of rational roots are all polynomial in $\size(f)$ (the sparse 
size of $f$). Furthermore, {\bf all} 
rational roots of $f$ can be computed within $\cO(\size(f)^{10})$ bit 
operations.\footnote{The exponent was not stated explicitly in 
\cite{lenstra} but, via \cite{lll}, can easily be derived from the description 
of the algorithm given there.} \qed 
\end{thm}  
\noindent 
Note that the complexity bound above does {\bf not} follow directly from 
the famous polynomial-time factoring algorithm of Lenstra, Lenstra, and 
Lovasz \cite{lll}: their result has complexity polynomial in the degree of 
$f$, as well as $\size(f)$. Also, Lenstra actually derived a more general 
version of the theorem above which applies to finding all bounded 
degree factors of a univariate polynomial over any fixed algebraic number 
field \cite{lenstra}. Interestingly, the analogue of theorem 
\ref{thm:lenstra} for the {\bf SLP size} is an open problem and, like theorem 
\ref{thm:plai} and corollary \ref{cor:realuni}, has considerable impact within 
complexity theory (see theorem \ref{thm:tau} of section \ref{sec:int} for the 
full statement). 

Curiously, there is currently no known analogue of 
theorem \ref{thm:lenstra} for {\bf systems} of multivariate polynomials. The 
main reason is that the most naive generalizations easily lead to various 
obstructions and even some unsolved problems in number theory.  
For example, as of mid-2000, it is still unknown whether deciding the 
existence of a rational root for $y^2\!=\!ax^3+bx+c$ is even Turing-decidable. 
Thus, the first obvious restriction to make, following the notation of the 
last two sections, is to consider only those $F$ where $Z_F$ is finite. But 
even then there are complications: 
\begin{itemize} 
\item[{\bf Q$_{\mathbf{1}}$}]{The number of integral roots of $F$ can actually 
be exponential in the sparse size of $F$: A simple example is the system 
\mbox{$(\prod^D_{i=1}(x_1-i),\ldots, \prod^D_{i=1}(x_n-i))$,} 
which has $D^n$ integral roots and a sparse size of $\cO(nD\log D)$. \qed }  
\item[{\bf Q$_{\mathbf{2}}$}]{For $n\!>\!1$, the 
integral roots of $F$ can have coordinates with bit-length exponential in 
$\size(F)$, thus ruling out one possible source $\np$ certificates: 
For example, the system $(x_1-2,x_2-x^2_1,\ldots,x_n-x^2_{n-1})$ 
has sparse size $\cO(n)$ but has $(1,2,\ldots,2^{2^{n-2}})$ as a root. \qed } 
\end{itemize}  
So it appears that restricting to deciding the existence of rational 
roots, instead of finding them, may be necessary for sub-exponential 
complexity. Nevertheless, these difficulties may disappear when $n$ is fixed: 
even the case $n\!=\!2$ is open.  

As for simple complexity upper bounds, the efficient deterministic algorithms 
of section \ref{sec:complex} can easily be converted to $\pspa$ 
algorithms for finding all rational points within the zero-dimensional part of 
an algebraic set. However, we will use a different approach to place this 
problem within an even lower complexity class: testing the densities of primes 
with certain properties. 

First note that averaging over many primes (as 
opposed to employing a single sufficiently large prime) is essentially 
unavoidable if one wants to use mod $p$ root counts to decide the existence 
of rational roots. For example, from basic quadratic residue theory \cite{hw}, 
we know that the number of roots $x^2_1+1$ mod $p$ is {\bf not} constant for 
sufficiently large prime $p$. 
Similarly, Galois-theoretic considerations are also necessary before using 
mod $p$ root counts to decide feasibility over $\Q$.
\setcounter{ex}{2}
\begin{ex}
\label{ex:cool}  
Take $m\!=\!n\!=\!1$ and $F\!=\!f_1\!=(x^2_1-2)(x^2_1-7)(x^2_1-14)$. Clearly, 
$F$ has no rational roots. However, it is easily checked via the 
Jacobi symbol \cite{hw,bs} that $F$ has a root 
mod $p$ for {\bf all} primes $p$. In particular, note that the Galois group 
here is not transitive: there is no automorphism of $\overline{\Q}$ which 
fixes $\Q$ and sends, say, $\sqrt{2}$ to $\sqrt{7}$. 
\end{ex}
 
So let us now state a precursor to our method for detecting rational roots:  
Recall that $\pi(x)$ denotes the number of primes $\leq\!x$. 
Let $\pi_F(x)$ be the variation on $\pi(x)$ where we instead count the number 
of primes $p\leq\!x$ such that the reduction of $F$ mod $p$ has a root in 
$\Z/p\Z$, and let $\#$ denote set cardinality. 
\begin{main}
\label{main:start} 
(See \cite[Thm.\ 2]{jcs}.) 
Following the notation of sections \ref{sec:complex} and \ref{sec:real}, 
assume now that the coefficients of $F$ are integers. Let $K$ be the 
field $\Q(x_i \; | \; (x_1,\ldots,x_n)\!\in\!Z_F \ , \ 
i\!\in\!\{1,\ldots,n\})$. Then the truth of GRH implies the two statements for 
all $x\!>\!33766$:
\begin{enumerate}
\item{ Suppose $\infty\!>\!\#Z_F\!\geq\!2$ and $\mathrm{Gal}(K/\Q)$ acts
transitively on $Z_F$. Then
\[ \frac{\pi_F(x)}{\pi(x)}< \left(1-\frac{1}{V_F}
\right)\left(1+\frac{(V_F!+1)\log^2 x + V_F!V_F\cO(V_F+\sigma(h_F))\log x}
{\sqrt{x}}\right).\] }
\item{ Suppose $\#Z_F\!\geq\!1$. Then independent of $\mathrm{Gal}(K/\Q)$,
we have
\[ \frac{\pi_F(x)}{\pi(x)}> \frac{1}{V_F}(1-b(F,x)),\] }
\end{enumerate}
where $0\!\leq\!b(F,x)\!<\!\frac{4V_F\log^2 x+
V^2_F\cO(V_F+\sigma(h_F)+nV_F\sigma(h_F)/\sqrt{x})\log
x}{\sqrt{x}}$ and $0\!\leq\!\sigma(h_F)\!=$\\
$\cO\left(M_F[\sigma(F)+m(n\log D+\log m)]+n^2V_F\log 
D\right)$. Better still, we have 
$\sigma(h_F)\!=\!\cO(M_F\sigma(F)+n^2V_F\log D)$ when $m\!\leq\!n$. 
\qed 
\end{main} 
\noindent 
The upper bound from assertion (1) appears to be new, and the lower bound from 
assertion (2) significantly improves earlier  
bounds appearing in \cite{hnam,morais,peter} which were polynomial in 
$D^n$. Explicit formulae for the above asymptotic estimates appear 
in \cite[Remarks 9 and 10]{jcs}. 

Theorem \ref{main:start} thus presents the first main 
difference between feasibility testing over $\C$ and $\Q$: from  
theorem \ref{main:koi}, we know that 
the mod $p$ reduction of $F$ has a root
in $\Z/p\Z$ for a density of primes $p$ which is either positive or zero, 
according as $F$ has a root in $\C$ or not. 
The corresponding gap between densities
happened to be large enough for Koiran's randomized 
oracle algorithm to decide feasibility over $\C$ (cf.\ 
section \ref{sec:complex}). (We point out that Koiran's 
algorithm actually relies on the behavior of the function  
$N_F$ defined below, which is more amenable than that of $\pi_F$.) 
On the other hand, assertion (1) of theorem \ref{main:start} tells 
us that the mod $p$ reduction of $F$ has a root in $\Z/p\Z$ for 
a density of primes which is either $1$ or $1-\frac{1}{V_F}$, 
according as $F$ has, or {\bf strongly} fails to have, a rational root.

Unfortunately, the convergence of $\frac{\pi_F(x)}{\pi(x)}$ to
its limit is unfortunately too slow to permit any obvious algorithm using
subexponential work. However, via a Galois-theoretic trick (cf.\ 
theorem \ref{thm:galois} below) we can nevertheless place rational root 
detection in a lower complexity class than previously known. 
\begin{main}
\label{main:riemann}
\cite{jcs}  
Following the notation and assumptions of theorem \ref{main:start}, 
assume further that $F$ fails to have 
a rational root $\Longleftrightarrow [Z_F\!=\!\emptyset$ or 
$\mathrm{Gal}(K/\Q)$ acts transitively on $Z_F]$. Then the truth of GRH 
implies that deciding whether $F$ has a rational root 
can be done in polynomial-time, given access to an oracle in $\np^\np$, i.e., 
within the complexity class $\crap$.  
Also, we can check the emptiness and  
finiteness of $Z_F$ unconditionally (resp.\ assuming GRH) within 
$\pspa$ (resp.\ $\am$). \qed 
\end{main}
\noindent 
The new oracle can be summarized as follows: Given any $F$ 
and a finite subset $S\!\subset\!\N$, our oracle instantaneously 
tells us whether or not there is a prime $p\!\in\!S$ such that 
the mod $p$ reduction of $F$ has {\bf no} roots in $\Z/p\Z$. 

Part of the importance of oracle-based algorithms, such as the 
one above or the algorithm from section \ref{sec:complex}, 
is that it could happen that $\pp\!\neq\!\np$ but the 
higher complexity classes we have been alluding to all collapse 
to the same level. For example, while it is known that 
$\mathbf{NP}\!\cup\!\mathbf{BPP}\!\subseteq\!\mathbf{AM}\! 
\subseteq\!\pp^{\np^\np}\!\subseteq\!
\np^{\np^\np}\!\subseteq\!\cdots\!\subseteq\! 
\mathbf{PSPACE}$, the properness of each inclusion is still unknown 
\cite{zachos,lab,arith,papa}. 

The algorithm for theorem \ref{main:riemann} is almost as simple 
as the algorithm for theorem \ref{main:koi} given earlier, and 
can be outlined as follows: 

\begin{itemize} 
\item[{\bf Step 0}]{ Let $N_F(x)$ denote the 
{\bf weighted} version of $\pi_F(x)$ where we instead sum the 
total number of roots in $\Z/p\Z$ of the mod $p$ reductions of $F$ 
over {\bf all} primes $p\!\leq\!x$.} 
\item[{\bf Step 1}]{ Let $t^*_0$ be an integer just large enough so that 
$t^*_0\!>\!33766$ and $b(F,t^*_0)\!<\!\frac{1}{10}$.}  
\item[{\bf Step 2}]{Estimate, via a constant-factor approximate counting 
algorithm of Stockmeyer \cite{stock}\footnote{
Stockmeyer's algorithm actually applies to any function from the 
complexity class $\#\pp$, and it is easily verified that 
$N_F$ and $\pi_F$ lie within this class.}, both $N_F(t^*_0)$ and 
$\pi_F(t^*_0)$ within a factor of $\frac{9}{8}$, using polynomially many calls 
to our $\np^\np$ oracle. Call these approximations 
$\bar{N}$ and $\bar{\pi}$ respectively. } 
\item[{\bf Step 3}]{ If $\bar{N}\!\leq\!(\frac{9}{8})^2\bar{\pi}$, declare 
$Z_F\cap\Q^n$ empty. Otherwise, declare $Z_F\cap\Q^n$ \mbox{nonempty. \qed} } 
\end{itemize} 
\noindent 
That our algorithm runs in polynomial time follows easily from our quantitative 
estimates from theorem \ref{main:start} and an analogous estimate 
for $N_F(x)$ (which also depends on GRH) from \cite{jcs}. The same holds for 
the correctness of our algorithm. 

Let us now close with some remarks on the strength of our last two theorems: 
First note that our restrictions on the input $F$ are actually rather gentle. 
In particular, if one assumes $m\!\geq\!n$ and fixes the 
monomial term structure of $F$, then it follows easily 
from the theory of resultants \cite{gkz94,introres,gcp} that, for a 
generic choice of the coefficients, $F$ will have only finitely many 
roots in $\Cn$. (See section \ref{sec:int} for our definition of 
generic.) Furthermore, it is quite frequently the case that our hypothesis 
involving $Z_F$ and $\mathrm{Gal}(K/\Q)$ holds when $F$ fails to have a
rational root.
\begin{thm} 
\label{thm:galois} 
\cite[Thm.\ 4]{jcs} 
Following the notation above, fix the monomial term structure of $F$ and 
assume further that $m\!\geq\!n$ and the coefficients of $F$ are integers of 
absolute value $\leq\!c$. Then the fraction of such $F$ with
$\mathrm{Gal}(K/\Q)$ acting transitively
on $Z_F$ is at least $1-\cO(\frac{\log c}{\sqrt{c}})$.
Furthermore, we can check whether $\mathrm{Gal}(K/\Q)$ acts
transitively on $Z_F$ within $\expt$ or, if one assumes GRH, within
$\pp^{\np^\np}$. \qed 
\end{thm} 
\noindent 
Thus, if the monomial term structure of $F$ is such that
$\#Z_F\!\neq\!1$ for a generic choice of the coefficients, it
easily follows that at least a fraction of $1-\cO(\frac{\log c}{\sqrt{c}})$ of
the $F$ specified above also have no rational roots. The case
where the monomial term structure of $F$ is such that $\#Z_F\!=\!1$
for a generic choice of the coefficients is evidently quite rare,
and will be addressed in future work.
\begin{rem}
A stronger version of the $m\!=\!n\!=\!1$ case of theorem \ref{thm:galois}
(sans complexity bounds) was derived by Gallagher in
\cite{gala}. The $m\!\geq\!n\!>\!1$ case follows from a combination of our
framework here, the Lenstra-Lenstra-Lovasz (LLL) algorithm \cite{lll}, and
an effective version of Hilbert's Irreducibility Theorem from \cite{cohen}. 
\qed 
\end{rem}

As we have seen, transferring conditional speed-ups from 
$\C$ to $\Q$ presents quite a few subtleties, and these are covered at length 
in \cite{jcs}. We also point out that there appears to be no obstruction to 
extending our algorithm above to detecting rational points over any 
fixed number field, within the same complexity bound. This will be pursued in 
future work. 

\subsection{Related Results Over $\Q$} 
We have mainly concentrated on the complexity of detecting rational points on 
certain zero-dimensional algebraic sets, which has been a somewhat overlooked 
topic. Indeed, while a $\pspa$ complexity bound for this problem could have 
been derived via, say, the techniques of \cite{chigo} no later than 1984, 
there appears to be no explicit 
statement of this fact. In any event, that a large portion of this 
problem can be done within the polynomial hierarchy appears to be new. 

On the other hand, for algebraic sets of positive dimension, even the 
decidability of feasibility over $\Q$ is open. That the study of rational 
points on higher-dimensional varieties has been, and continues to be, 
intensely studied by some of the best number theorists and algebraic geometers 
is a testament to the difficulty of this problem. Current work on finding 
rational points has thus focused on characterizing (in terms of the underlying 
complex geometry) when a variety has infinitely many rational points, and how 
and where density of rational points can appear. 

For example, it was unproved until the work of Faltings in 1983 
\cite{faltings,bomb} that algebraic curves of genus\footnote{We will 
use {\bf geometric} (as opposed to arithmetic) genus throughout 
this paper. These definitions can be found in \cite{hart,miranda}.} 
$\geq\!2$ have only finitely many rational points. (This fact was originally 
conjectured by L.\ J.\ Mordell in 1922.) The seminal work of Lang and 
Vojta has since lead to even deeper connections between the distribution of 
rational points and the geometry of the underlying complex manifold 
\cite{vojta,lang}. More recently, highly refined quantitative results 
(some depending on conjectures of Lang) on the density of rational points on 
certain varieties have appeared (see, e.g., \cite{manin,pacelli,tschinkel} and 
the references therein). 

This is of course but a fragment of the wealth of current 
active research on rational points, and we have yet to 
speak of the complexity of finding integral points. 

\section{Effective Siegel Versus Detecting Integral Points on Surfaces} 
\label{sec:int} 
The final results we present regard the computational complexity of 
certain problems involving integral points on varieties of dimension 
$\geq\!1$. We will strike a path leading to a relation between 
height bounds for integral points on algebraic plane curves and certain 
Diophantine prefixes in $\leq\!4$ variables, e.g., sentences of the 
form \[ \exists u\!\in\!\N \; \forall x\!\in\!\N \; \exists y\!\in\!\N \; 
f(u,x,y)\!\stackrel{?}{=}0.\] (The last sentence is an example of the prefix 
$\exists\forall\exists$, and we will casually refer to various quantified sentences in 
this way.) We then conclude with some evidence for the 
undecidability of Hilbert's Tenth Problem in three variables (theorem 
\ref{main:h3}). 

We first note that Diophantine complexity has quite a rich theory 
already in one variable. 
\begin{thm} 
\label{thm:tau}
\cite[Thm.\ 3, pg.\ 127]{bcss} 
Let $\tau(f)$ denote the SLP size of $f\!\in\!\Z[t]$, starting from 
the sequence $\{1,t,\ldots\}$. Suppose there exists an absolute 
constant $C_2\!>\!0$ such that for all $f$, the number of 
integral roots of $f$ is bounded above by $(\tau(f)+1)^{C_2}$. 
Then $\pp_\C\!\neq\!\np_\C$.\footnote{i.e., the analogue of the 
$\pp\!\neq\!\np$ conjecture for the BSS model over $\C$ would 
be settled.} \qed 
\end{thm} 
\noindent  
In short, a deeper understanding 
of the SLP encoding (cf.\ section \ref{sec:real}) over $\Z$ would have a 
tremendous impact in complexity theory. 

Via the sparse encoding, the study of integral roots for polynomials in two 
variables leads us to similar connections with important complexity classes.  
\begin{thm}
\label{thm:biint}
\cite{adleman} 
Deciding whether $ax^2+by\!=\!c$ has a root $(x,y)\!\in\!\N^2$, for an 
arbitrary input $(a,b,c)\!\in\!\N^3$, is $\np$-complete relative 
to the sparse encoding. i.e., there is an algorithm for this 
problem with bit complexity polynomial in $\log(abc)$ iff 
$\pp\!=\!\np$. \qed
\end{thm}
\noindent 
Note that we hit the class $\np$ rather quickly: quadratic polynomials  
(or genus zero curves)\footnote{It will be convenient to describe bivariate 
polynomials in terms of their underlying complex geometry, and we will 
do so freely in this section.} are enough. The case of higher degree 
polynomials is much less understood. To see this, let us denote the following 
problem by $\htp(n)$:\\ 
\vspace{-.4cm}
\begin{center}
``Decide whether an arbitrary
$f\!\in\!\Z[x_1,\ldots,x_n]$ has a root in $\N^n$ or not.''\footnote{
Hilbert's Tenth Problem in $n$ variables is actually the simplification of 
$\htp(n)$ where we seek roots in $\Zn$. However, for technical reasons, it is 
more convenient to deal with $\htp(n)$. }
\end{center}
\noindent  
(So our last theorem can be rephrased as the $\np$-hardness of $\htp(2)$  
for quadratic polynomials.) It is then rather surprising that as of mid-2000, 
the decidability of $\htp(2)$ is still open, even for general polynomials of 
degree $4$ (or general curves of genus $2$). 

Alan Baker has conjectured 
\cite[Section 5]{jones81} that the analogue $\htp(2)$ for 
$\Z^2$ is decidable. 
More concretely, the decidability of $\htp(2)$ is known in certain special 
cases, and these form a significant part of the applications of Diophantine 
approximation and arithmetic geometry. To describe the 
known cases, it is convenient to introduce the following functions. 
\begin{dfn}
\label{dfn:big}
Following the notation of sections \ref{sec:complex} and 
\ref{sec:real}, define the function $\biggy_\N : \Z[x_1,x_2] \longrightarrow 
\N\cup\{0,\infty\}$ by letting $\biggy_\N(f)$ be the
supremum of $\max\{|r_1|,|r_2|\}$ as $(r_1,r_2)$
ranges over $\{(0,0)\}\cup (Z_f\cap\N^2)$. The function $\biggy_\Z(f)$ is 
defined similarly, simply letting $(r_1,r_2)$ range over $\{(0,0)\}\cup 
(Z_f\cap\Z^2)$ instead. \qed 
\end{dfn} 
Parallel to $\htp(n)$ and its analogue over $\Zn$, 
the computability of $\biggy_\N$ implies the computability of 
$\biggy_\Z$. (Simply consider the substitution $f(x,y) \mapsto 
f(-x,-y)f(-x,y)f(x,-y)f(x,y)$.) The other direction is actually 
not trivial: there is nothing stopping a curve from having 
infinitely many integral points {\bf outside} of the 
first quadrant, thus obstructing any useful bound for $\biggy_\Z$ 
from being a useful bound for $\biggy_\N$. 

The computability of $\biggy_\N$ would of course imply the 
decidability of $\htp(2)$. However, as of mid-2000, even 
the computability of $\biggy_\Z$ is, with a few exceptions, known only 
for those $f$ where $Z_f$ falls into one of the following cases: certain genus 
zero curves \cite{poulaki}, all genus one curves \cite{bakercoates}, certain 
genus two curves \cite{grant,poon}, Thue curves \cite{bakert}, and curves in 
super-elliptic form \cite{bakerh,brindza}. (These also happen to be the 
only cases for which the decidability of $\htp(2)$ is known.) 
For example, it is known that for any polynomial equation of the form
\[ y^2=a_0+a_1x+a_2x^2+a_3x^3, \]
where $a_0,a_1,a_2,a_3\!\in\!\Z$ and $a_0+a_1x+a_2x^2+a_3x^3$ has three
distinct complex roots, all integral solutions must satisfy
\[|x|,|y|\leq \exp((10^6c)^{10^6}),\]
where $c$ is any upper bound on $|a_0|,|a_1|,|a_2|,|a_3|$ \cite{bakertran}. 
(More recent improvements of this bound can be found in \cite{schmidt}.) 
\begin{rem} 
An interesting related conjecture of Steve Smale \cite{steve} is 
that if a plane curve of positive genus has an integral point, then it  
must have an integral point of height singly exponential in the dense size of 
the defining polynomial. (See below for the definition of 
dense size.) \qed  
\end{rem} 

Of course, one may still worry whether $\biggy_\Z$ can be computable 
without $\biggy_\N$ being computable. We can resolve this 
as follows: 
\begin{thm} 
\label{main:equiv} 
The function $\biggy_\N$ is computable $\Longleftrightarrow \biggy_\Z$ is 
computable. 
\end{thm} 
\noindent 
The proof follows easily from theorem \ref{thm:silvb} of the next section, 
which describes the distribution of integral points within the real part 
of a complex curve. In spite of theorem \ref{main:equiv}, it is still unknown 
whether replacing $\Z^2$ by $\N^2$ makes a significant difference in the 
complexity of $\htp(2)$. However, via theorem \ref{thm:siegel} of the 
next section, we can prove that a variant of $\htp(2)$ is closely 
related to detecting infinitudes of integral points on plane curves.
\begin{thm}
\label{main:25}
Let $\ratcurve(3)$ denote the problem of deciding whether a 
(geometrically irreducible, possibly singular) 
genus zero curve in $\C^3$ defined over $\Z$ contains a 
point in $\N^3$. Also let $\htp^\infty(2)$ denote the 
problem of deciding whether an arbitrary $f\!\in\!\Z[x_1,x_2]$ has 
infinitely many roots in $\N^2$. Then $\ratcurve(3)$ decidable 
$\Longrightarrow \htp^\infty(2)$ decidable. 
\end{thm} 
\noindent 
We note that the input for $\ratcurve(3)$ is given as usual: 
a set of polynomials in $\Z[x_1,x_2,x_3]$ defining the curve in 
question. 
Curiously, the decidability of $\ratcurve(3)$, $\htp^\infty(2)$, and 
their analogues over $\Z$ are all unknown, in spite of Siegel's Theorem. 
(Siegel's Theorem \cite{siegel} is a famous result from 1934 partially 
classifying those curves with infinitely many integral points.) A  
refined version of Siegel's Theorem appears as theorem \ref{thm:siegel} of 
the next section.

The preceding results can all be considered as variations on the study of 
the Diophantine prefixes $\exists$ and $\exists\exists$. So to prove more 
decisive results it is natural to study subtler combinations of quantifiers. 
In particular, we will show that the prefix 
$\exists\forall\exists$ can be solved (almost always) within the polynomial 
hierarchy. To make this more precise, let us make two quantitative 
definitions: 
When we say that a statement involving a set of parameters 
$\{c_1,\ldots,c_N\}$ is true {\bf generically}\footnote{ 
We can in fact assert a much stronger condition, but this one suffices for 
our present purposes.}, we will 
mean that for any $M\!\in\!\N$, the statement fails for at most 
$\cO(N(2M+1)^{N-1})$ of the $(c_1,\ldots,c_N)$ lying in $\{-M,\ldots,M\}^N$. 
Also, for an algorithm with a polynomial $f\!\in\!\Z[v,x,y]$ as input, 
speaking of the {\bf dense encoding} will simply mean 
measuring the input size as $D+\sigma(f)$, where $D$ (resp.\ $\sigma(f)$) is 
the total degree (resp.\ maximum bit-length of a coefficient) of $f$.  
\begin{thm}
\cite{jcs}
\label{main:pepper} 
Fix the Newton polytope $Q$ of a polynomial $f\!\in\!\Z[v,x,y]$  
and suppose that $Q$ has 
at least one integral point in its interior.\footnote{So, among other things, 
we are assuming $Q$ is $3$-dimensional.} Assume further that we measure 
input size via the dense encoding. Then, for a generic choice of 
coefficients depending only on $Q$, we can decide whether 
$\exists v \; \forall x \; \exists y \; f(v,x,y)
\!=\!0$ 
(with all three quantifiers ranging over $\N$ or $\Z$) within $\conp$. 
Furthermore, we can check whether an input $f$ has generic coefficients within 
$\nc$. \qed 
\end{thm} 
\noindent
The hierarchy of complexity classes $\nc$ simply consists of those problems 
in $\pp$ which admit efficient parallel algorithms (see \cite{papa} for a full 
statement). Roughly speaking, 
deciding the prefix $\exists\forall\exists$ is equivalent to determining 
whether an algebraic surface has a slice (parallel to the 
$(x,y)$-plane) densely peppered with integral points, and we have thus 
shown that this problem is tractable for most inputs.  Whether 
$\conp$-completeness persists relative to the {\bf sparse} encoding 
remains an open question. 

It is interesting to note that the exceptional case to our 
algorithm for $\exists\forall\exists$ judiciously contains an extremely hard 
number-theoretic problem: the prefix $\exists\exists$ or, equivalently, 
$\htp(2)$. (That $\Z[v,y]$ lies in our exceptional locus is easily checked.) 
More to the point, James P.\ Jones has conjectured \cite{jones81} that the 
decidabilities of the prefixes 
$\exists\forall\exists$ and $\exists\exists$, quantified over $\N$, are 
equivalent. Thus, while we have not 
settled Jones' conjecture, we have at least shown that the decidability of 
$\exists\forall\exists$ now hinges on a sub-problem much closer to 
$\exists\exists$. 
 
Call an algebraic surface $Z\!\subset\!\C^4$ {\bf specially ruled} 
iff it is a bundle of genus zero curves fibered over a genus zero 
curve in the $(u,v)$-plane (coordinatizating $\C^4$ by $(u,v,x,y)$). 
The proof of theorem \ref{main:pepper} is primarily based on a geometric trick 
which easily extends to the prefix $\exists\exists\forall\exists$. In 
particular, we also have the following result. 
\begin{thm}
\label{main:h3} 
At least one of the following two statements is {\bf false}:
\begin{enumerate}
\item{The function $\biggy_\N$ is Turing-computable.}
\item{The Diophantine sentence
\[ \exists u\!\in\!\N \ \ \exists v\!\in\!\N \ \ \forall
x\!\in\!\N \ \
\exists y\!\in\!\N  \ \ f(u,v,x,y)\!\stackrel{?}{=}\!0 \]
is decidable in the special case where the underlying
$3$-fold $Z_f$ contains a specially ruled surface. }
\end{enumerate}
In particular, $\htp(3)$ is a special case of the problem
mentioned in statement (2). 
\end{thm}
\noindent 
A slightly stronger version of theorem \ref{main:h3} appears in 
\cite{tcs} and, for the convenience of the reader, we supply a more 
streamlined proof in section \ref{sub:h3}. 
We thus now have (applying theorem \ref{main:equiv}) a weak version of the 
following implication: $\biggy_\Z$ computable $\Longrightarrow \htp(3)$ 
undecidable. 

Since Matiyasevich and Robinson have shown 
that $\exists\exists\forall\exists$ is undecidable (when 
all quantifiers range over $\N$) \cite{matrob}, our last theorem can also 
be interpreted as a restriction of this undecidability to a particular 
subset of the general problem. Whether this subproblem can 
be completely reduced to $\htp(3)$ is therefore of the utmost interest. 

\subsection{Related Work Over $\N$ and $\Z$} 
\label{sub:relatedint}
We first point out that the decidability of $\exists\forall\exists$ was an open 
problem and, in spite of theorem \ref{main:pepper}, remains open for {\bf 
arbitrary} inputs. We also note that our algorithm for (most of) 
$\exists\forall\exists$ is based on an important result of Tung for the prefix 
$\forall\exists$.
\begin{tungy}
\cite{tungcomplex} 
Deciding the quantifier prefix $\forall\exists$ (with all quantifiers 
ranging over $\N$ or $\Z$) is $\conp$-complete relative to the dense 
encoding. \qed 
\end{tungy} 
\noindent 
The decidability of $\forall\exists$ (over $\N$ and $\Z$) was first derived by 
James P.\ Jones in 1981 \cite{jones81}. 
The algorithms for $\forall\exists$ alluded to in Tung's Theorem are based 
on some very elegant algebraic facts due to Jones, Schinzel, and Tung. 
We illustrate one such fact for the case of $\forall\exists$ over $\N$. 
\begin{jst}
\cite{jones81,schinzel,tungcomplex}
Given any $f\!\in\!\Z[x,y]$, we have that
$\forall x\; \exists y\; f(x,y)\!=\!0$
iff all three of the following conditions hold:
\begin{enumerate}
\item{The polynomial $f$ factors into the form
$f_0(x,y)\prod^j_{i=1}(y-f_i(x))$ where $f_0(x,y)\!\in\!\Q[x,y]$
has {\bf no} zeroes in the ring $\Q[x]$, and for all $i$,
$f_i\!\in\!\Q[x]$ and the leading coefficient of $f_i$ is positive.}
\item{$\forall x\!\in\!\{1,\ldots,x_0\} \; \exists
y\!\in\!\N$ such that $f(x,y)\!=\!0$, where $x_0\!=\!\max\{s_1,\ldots,s_j\}$,
and for all $i$, $s_i$ is the sum of the squares of the coefficients
of $f_i$.}
\item{Let $\alpha$ be the least positive integer such that
$\alpha f_1,\ldots,\alpha f_j\!\in\!\Z[x]$ and set 
$g_i\!:=\!\alpha f_i$ for all $i$. 
Then the {\bf union} of the solutions of the following $j$ congruences
$g_1(x)\!\equiv\!0 \ (\mod \ \alpha), \ldots, g_j(x)\!\equiv\!0 \ (\mod \ 
\alpha)$ is {\bf all} of $\Z/\alpha\Z$. \qed }
\end{enumerate}
\end{jst}
\noindent 
The analogue of the JST Theorem over $\Z$ is essentially the same, save for the 
absence of condition (2), and the removal of the sign check in condition (1) 
\cite{tungcomplex}. 

The study of the decidability of Diophantine prefixes dates back to 
\cite{oldmat,matrob,jones81}, and \cite{hilbert10,tungnew,stoc99,jcs} give 
some of the most recent results. Of course, as we have seen above, there is 
still much left to be done, and we hope that this paper sparks the 
interests of other researchers. 

In particular, the precise complexity of checking whether 
an input $f\!\in\!\Z[u,v,x,y]$ satisfies the hypothesis of statement 
(2) of theorem \ref{main:h3} is unknown. (The decidability of 
this problem is at least known, and there are more restricted versions of 
(2) which can be checked within $\nc$ \cite{tcs}.) 
The author conjectures that this hypothesis can in fact be decided within 
$\nc$, relative to the dense encoding.

More to the point, it is curious that the complexity of 
deciding whether a given curve has infinitely many integral 
points is also open. 
The best result along these lines 
is the following refined version of Siegel's Theorem: 
\begin{thm}
\label{thm:siegel}
\cite{wow} 
Following the notation of sections \ref{sec:complex} and \ref{sec:real}, 
suppose $f\!\in\!\Z[x_1,x_2]$ is such that $Z_f$ is a geometrically 
irreducible curve. Then $Z_f\cap\Z^2$ is infinite $\Longleftrightarrow$ all of 
the following three conditions are satisfied: 
\begin{enumerate}
\item[(a)]{ $Z_f$ has genus $0$,} 
\item[(b)]{ $Z_f\cap\Z^2$ contains at least one non-singular point, and}  
\item[(c)]{ the highest degree part of $f$ has either (i) exactly one root 
in $\Pro^1_\C$ (necessarily rational) or (ii) has exactly two distinct roots 
in $\Pro^1_\C$ {\bf and} they are both real. \qed }
\end{enumerate} 
\end{thm} 
\noindent 
Joseph H.\ Silverman has pointed out that this result may 
already be known to experts in algebraic curves. 
Another curious fact regarding Siegel's theorem is that it still has no proof 
which settles the computability of $\biggy_\Z$. 

A useful result arising from Silverman's proof of theorem 
\ref{thm:siegel} is the following solution to a conjecture 
of the author from \cite{tcs}: 
\begin{thm} 
\cite{wow}
\label{thm:silvb}
Let $W$ be any geometrically irreducible curve in $\C^2$ defined over 
$\Z$ possessing infinitely many integral points. Let $W'$ be any unbounded 
subset of $W\cap\R^2$. Then $W'$ contains infinitely many integral points. 
\qed 
\end{thm} 
\noindent 
This result, combined with a little computational algebraic geometry, 
provides the proof of theorem \ref{main:equiv} and the details appear 
in section \ref{sub:proofint}. 

As for more general relations between $\htp(n)$ and its analogue over 
$\Zn$, it is easy to see that the decidability of $\htp(n)$ implies 
the decidability of its analogue over $\Zn$. Unfortunately, 
the converse is currently unknown. Via Lagrange's Theorem 
(that any positive integer can be written as a sum of four 
squares) one can easily show that the {\bf un}decidability of 
$\htp(n)$ implies the {\bf un}decidability of the analogue 
of $\htp(4n)$ over $\Zn$. More recently, Zhi-Wei Sun has 
shown that the $4n$ can be replaced by $2n+2$ \cite{sun}.  

\section{Proofs of Our Main Technical Results} 
\label{sec:proofs} 
For the convenience of the reader, let us briefly distinguish what is new 
and/or recent: To the best of the author's knowledge, theorems 
\ref{main:complex}, \ref{main:koi}, \ref{main:height}, \ref{main:size}, 
\ref{main:unired}, \ref{main:equiv}, and \ref{main:25}, and corollary 
\ref{cor:realuni} have not appeared in print before. Also, although theorem 
\ref{main:equiv} was conjectured, along with a plan of attack, in \cite{tcs}, 
its full proof has not appeared before. Finally, while preliminary versions 
of theorems \ref{main:height} \and \ref{main:unired} appeared earlier 
in \cite{gcp}, their corresponding height bounds are new. 

As for the remaining results, they have either already appeared, or are about 
to appear, in the references listed in their respective statements.  

Our proofs will thus focus on results over our ``outlying'' rings: 
$\C$ and $\Z$. 

\subsection{Proofs of Our Results Over $\C$: Theorems \ref{main:complex}, 
\ref{main:height}, \ref{main:size}, \ref{main:unired}, and 
\ref{main:koi}}\mbox{}\\
\label{sub:proofcomplex}
While our proof of theorem \ref{main:koi} will not directly 
require knowledge of resultants, our proofs of theorems \ref{main:complex}, 
\ref{main:height}, \ref{main:size}, and \ref{main:unired} are based on the 
{\bf toric resultant}.\footnote{Other 
commonly used prefixes for 
this modern generalization of the classical resultant \cite{vdv} include:
sparse, mixed, sparse mixed, $\cA$-, $(\cA_1,\ldots,\cA_k)$-, and Newton. 
Resultants actually date back to work Cayley and Sylvester in the 
19$^\thth$ century, but the toric resultant incorporates some 
combinatorial advances from the late 20$^\thth$ century. } 
This operator allows us to reduce all the computational algebraic geometry we 
will encounter to matrix and univariate polynomial arithmetic, with almost no 
commutative algebra machinery. We supply a precis on the toric resultant in 
the following section. 

As mentioned earlier, we will reduce the description of $Z_F$ to 
univariate polynomial factorization. Another trick we will use is to reduce 
most of our questions to finding isolated roots of polynomial systems where the 
numbers of equations and variables is the same. 

These geometric constructions are useful for the proof of theorem 
\ref{main:koi} 
as well, but more in a theoretical sense than in an algorithmic sense. As we 
will see in section \ref{sub:koi}, it is number theory which allows us to 
enter a lower complexity class, and univariate reduction is needed only for 
quantitative estimates. 

\subsubsection{Background on Toric Resultants}\mbox{}\\
\label{sub:res} 
Since we do not have the space to give a full introduction to 
resultants we refer the reader to \cite{emiphd,gkz94,introres} 
for further background. The necessary facts we need are all 
summarized below. In what follows, we let $[j]\!:=\!\{1,\ldots,j\}$. 

Recall that the {\bf support}, $\pmb{\supp(f)}$, of a polynomial 
$f\!\in\!\C[x_1,\ldots,x_n]$ is simply the set of exponent vectors of the 
monomial terms appearing\footnote{We of course fix an ordering on the 
coordinates of the exponents which is compatible with the usual ordering 
of $x_1,\ldots,x_n$.} in $f$. The support of the 
{\bf polynomial system} $F\!=\!(f_1,\ldots,f_m)$ is simply the $m$-tuple 
$\pmb{\supp(F)}\!:=\!(\supp(f_1),\ldots,\supp(f_m))$. Let 
$\bar{\cA}\!=\!(\cA_1,\ldots,\cA_{m+1})$ be any $(m+1)$-tuple of non-empty 
finite subsets of $\Zn$ and set $\cA\!:=\!(\cA_1,\ldots,\cA_m)$. If we say 
that $F$ has {\bf support contained in} $\cA$ then we simply mean that 
$\supp(f_i)\!\subseteq\!\cA_i$ for all $i\!\in\![m]$.  
\begin{dfn}
\label{dfn:res}  
Following the preceding notation, suppose we can find line 
segments $[v_1,w_1],\ldots,[v_{m+1},w_{m+1}]$  
with $\{v_i,w_i\}\!\subseteq\!\cA_i$ for all $i$ and $\vol_m(L)\!>\!0$, where  
$L$ is the convex hull of $\{\bO,w_1-v_1,\ldots,w_{m+1}-v_{m+1}\}$. Then we 
can associate to $\bar{\cA}$ a unique (up to sign) irreducible polynomial 
$\res_{\bar{\cA}}\!\in\!\Z[c_{i,a} \; | \; i\!\in\![m+1] 
 \ , \ a\!\in\!\cA_i]$ with the following property: If we identify 
$\bar{\cC}\!:=\!(c_{i,a} \; | \; i\!\in\![m+1] \ , \ a\!\in\!\cA_i)$ 
with the vector of coefficients of a polynomial system $\bar{F}$ with support 
contained in $\bar{\cA}$ (and constant coefficients), then $\bar{F}$ has a 
root in $\Csn \Longrightarrow \res_{\bar{\cA}}(\bar{\cC})\!=\!0$. Furthermore, 
for all $i$, the degree of $\res_{\bar{\cA}}$ with respect to the coefficients 
of $f_i$ is no greater than $V_F$. \qed 
\end{dfn} 
\noindent 
We 
emphasize that the implication above does {\bf not} go both ways: 
the correct converse involves toric varieties \cite{gkz94,jpaa,gcp}.  
A consequence of the above definition is that the toric resultant 
applies mainly to systems of $n+1$ polynomials in $n$ 
variables. However, via a trick from the next section, this 
will cause no significant difficulties when we consider $m$ polynomials 
in $n$ variables. 

That the toric resultant can actually be defined as above is covered 
in detail in \cite{combiresult,gkz94}. 
There is in fact an exact formula for the degree of $\res$ 
with respect to the coefficients of $f_i$ involving {\bf mixed} volumes  
\cite{combiresult,gkz94}. Our simplified upper bound follow easily 
from the fact that mixed volume never decreases when the input polytopes 
are grown \cite{buza}. 

Another operator much closer to our purposes is the {\bf toric 
perturbation} of $F$. 
\begin{dfn}
\label{dfn:pert}  
Following the notation of definition \ref{dfn:res}, assume further that 
$m\!=\!n$, $\supp(F)\!=\!\cA$, and $\supp(F^*)\!\subseteq\!\cA$. 
We then define \[ \pert_{(F^*,\cA_{n+1})}(u)\!\in\!\C[u_a \; | \; 
a\!\in\!\cA_{n+1}]\] to be the coefficient of the term of 
\[ \res_{\bar{\cA}}(f_1-sf^*_1,\ldots,f_n-sf^*_n,\!\!\!\!\sum_{a\in 
\cA_{n+1}}\!\!\!\!\!u_ax_a)\!\in\!\C[s][u_a \; | 
\; a\!\in\!\cA_{n+1}]\] of {\bf lowest} degree in $s$. \qed 
\end{dfn} 
\noindent 
The constant term of the last resultant is a generalization of the 
classical {\bf Chow form} of a zero-dimensional variety \cite{vdv}. 
The consideration of the higher order coefficients is necessary 
when $Z_F$ is positive-dimensional. In particular, 
the geometric significance of $\pert$ can be summarized as follows: 
For a suitable choice of $F^*$, $\cA_{n+1}$, and $\{u_a\}$, 
$\pert$ satisfies all the properties of the polynomial $h_F$ from theorem 
\ref{main:height} in the special case $m\!=\!n$. In essence, $\pert$ is an 
algebraic deformation which allows us to replace the positive-dimensional 
part of $Z_F$ by a finite subset which is much easier to handle.
 
To prove theorems \ref{main:complex}, \ref{main:height}, and \ref{main:unired} 
we will thus need a good complexity  estimate for computing $\res$ and $\pert$. 
\begin{lemma}
\label{lemma:respert} 
Following the notation above, 
let $\cR_F$ (resp.\ $\cP_F$) be the number of deterministic arithmetic 
operations needed to evaluate 
$\res_{\bar{\cA}}$ (resp.\ $\pert_{(F^*,\cA_{n+1})}$) at any point in 
$\C^{k+n+1}$ (resp.\ $\C^{2k+n+1}$), where $\cA\!\subseteq\!\supp(F)$ and 
$\cA_{n+1}\!:=\!\{\bO,
e_1,\ldots,e_n\}$. Also let $r_F$ be the total degree of 
$\res_{\bar{\cA}}$ as a polynomial in the coefficients of $\bar{F}$  
Then $r_F\!\leq\!(n+1)V_F$, 
$\cR_F\!\leq\!(n+1)r_F\cO(M^{2.376}_F)$, and 
$\cP_F\!\leq\!(r_F+1)\cR_F+r_F(1+\frac{3}{2}\log r_F)$. 
Furthermore, $k\!\leq\!m(V_F+n)$ and 
$M_F\!\leq\!e^{1/8}\frac{e^n}{\sqrt{n+1}}V_F+\prod^n_{i=1} 
(p_i+2)-\prod^n_{i=1} (p_i+1)$, where $p_i$ is the length of the 
projection of $nQ_F$ onto the $x_i$-axis. 
(Note that $e^{1/8}\!\approx\!1.3315$.) 
\qed 
\end{lemma}
\noindent 
{\bf Proof:} 
The bound on $\cR_F$ (resp.\ $\cP_F$) follows directly from \cite{emican} 
(resp.\ \cite{gcp}), as well as a basic complexity result on the 
{\bf inverse discrete Fourier transform} \cite[pg.\ 12]{binipan}. 

The bound on $k$ follows by noting that $k\!\leq\!m\ell_F$, where 
$\ell_F$ is the number of lattice points in the polytope $Q_F$. 
By a classical lattice point count of Blichfeldt \cite{blich}, 
we obtain $\ell_F\!\leq\!V_F+n$ and we are done. 

As for the bound on $M_F$, we will observe a bit later that 
$M_F$ can be bounded above by the number of lattice points in 
the {\bf Minkowski sum}\footnote{The Minkowksi sum of any two 
subsets $A,B\subseteq\Rn$ is simply the set $\{a+b \; | \; 
a\in A \ , \ b\in B\}$.} $Q'_F\!:=\!nQ_F+\conv\{\bO,e_1,\ldots,e_n\}$. (This 
polytope is clearly contained in the polytope $(n+1)Q_F$ mentioned in 
theorem \ref{main:complex}.) Noting that 
$\frac{(n+1)^n}{n!}\!\leq\!e^{1/8}\frac{e^n}{\sqrt{n+1}}$ via 
Stirling's estimate \cite[pg.\ 200]{rudin}, 
and that the length of the projection of $Q'_F$ onto the $x_i$-axis is exactly 
$p_i+1$, our bound on $M_F$ follows immediately from another 
simple lattice point count \cite[Formula 3.11]{gw}. \qed 
\begin{rem} 
That $M_F\!=\!\cO(V_F)$ for fixed $n$ is immediate from our last lemma.  
Note also that $Q'_F$ is contained in the standard $n$-simplex 
scaled by a factor $nD+1$. Calling the latter polytope $\cQ_F$, 
it is clear that the number of lattice points in $\cQ_F$ 
is yet another upper bound on $M_F$. The latter lattice point count 
in turn has a simple explicit formula in terms of the 
binomial coefficient, and this is how we derived the 
crude bound on $M_F$ mentioned in section \ref{sub:relatedcomplex}. 
\qed 
\end{rem} 

Admittedly, such complexity estimates seem rather mysterious without  
any knowledge of how $\res$ and $\pert$ are computed. So let us 
now give a brief summary: 
The key fact to observe is that, in the best circumstances, one can express 
$\res$ as the determinant of a (square) sparse structured matrix 
$\cM_{\bar{\cA}}$ --- a {\bf toric resultant matrix} --- whose entries are 
either $0$ or polynomials 
in the coefficients of $\bar{F}$. (In fact, the $\cM_{\bar{\cA}}$  
we use will have every row equal to a permutation of the 
vector $v=(\cC_i,0,\ldots,0)$, where $\cC_i$ is the vector of coefficients 
of $f_i$ and $i$ (and the permutation) depends on the row.) These matrices 
have their origin in the study of certain spectral sequences \cite{gkz94} and 
there are now down-to-earth combinatorial algorithms for finding them 
\cite{emican,emiphd,emipanmat,emimoumat}. 
 
So the quantity $M_F$ in our theorems is nothing more than 
the number of rows (or columns) of $\cM_{\bar{\cA}}$. The bound 
on $M_F$ from our last theorem thus arises simply by applying the main 
algorithm from \cite{emican}, and observing that this particular 
construction of $\cM_{\bar{\cA}}$ creates a matrix row for every lattice point 
in a translate of the polytope $\conv(\cA_1+\cdots+\cA_{n+1})$. 
In particular, it is also the case that the deterministic arithmetic 
complexity of constructing $\cM_{\bar{\cA}}$ 
is dominated by $\cO(M_F\log n + n^2)$ \cite{new}, so we can 
henceforth ignore this construction in our complexity bounds. 
Better still, the quantity $M_F$ can be expected to admit even sharper 
upper bounds, once better algorithms for building toric resultant 
matrices are found. 

However, it is more frequently the case that $\res$ is but a {\bf divisor} of 
$\det \cM_{\bar{\cA}}$, and further work must be done. Fortunately, 
in \cite{emican,emiphd}, there are general randomized and deterministic 
algorithms for extracting $\res$. These algorithms work via subtle 
refinements of the classical technique of recovering the coefficients of a 
polynomial $g$ of degree $D$ by evaluating $g$ at $D+1$ points and then 
solving for the coefficients via a structured linear system. This accounts for 
the appearance of the famous linear algebra complexity exponent 
($\omega\!<\!2.376$), or simple functions thereof, in our 
complexity estimates. 

\subsubsection{The Proof of Theorem \ref{main:complex}}\mbox{}\\
Our algorithm can be stated briefly as follows: 
\begin{itemize} 
\item[{\bf Step 0}]{If $f_i$ is indentically $0$ for all $i$, 
declare that $Z_F$ has dimension $n$ and stop. Otherwise, 
let $i\!:=\!n-1$. } 
\item[{\bf Step 1}]{For each $j\!\in\![2k+1]$, 
compute an $(i+1)n$-tuple of integers\\ $(\eps_1(j),\ldots,\eps_n(j),
\eps_{(1,1)}(j),\ldots,\eps_{(i,n)}(j))$ via lemma \ref{lemma:probe} and 
the polynomial system (\ref{eq:probe}) below.} 
\item[{\bf Step 2}]{Via theorem \ref{main:height}, 
check if the polynomial system  
\begin{eqnarray} 
\label{eq:probe}  
\eps_1(j)f_1+\cdots+\eps_1(j)^mf_m+\eps_1(j)^{m+1}l_1+\cdots+\eps_1(j)^{m+i} 
l_i & = & 0 \notag \\ 
& \vdots  &  \\ 
\eps_n(j)f_1+\cdots+\eps_n(j)^mf_m+\eps_n(j)^{m+1}l_1+\cdots+
\eps_n(j)^{m+i} l_i & = & 0 \notag 
\end{eqnarray} 
has a root for more than half of the $j\!\in\![2k+1]$, 
where\\ 
$l_t\!:=\!\eps_{(t,1)}x_1+\cdots+\eps_{(t,n)}x_n$. } 
\item[{\bf Step 3}]{If so, declare that $Z_F$ has dimension $i$ and stop. 
Otherwise, if $i\!\geq\!1$, set $i\mapsto i-1$ and go to Step 1.} 
\item[{\bf Step 4}]{ Via theorem \ref{main:unired} and a univariate gcd 
computation, check if the system (\ref{eq:probe}) has a root which is 
also a root of $F$.} 
\item[{\bf Step 5}]{ If so, declare that $Z_F$ has dimension $0$ 
and stop. Otherwise, 
declare $Z_F$ empty and stop.} 
\end{itemize} 
\noindent  
Before analyzing the correctness of our algorithm, let us briefly 
clarify Steps 2 and 4. First let $G_{(j)}$ denote the polynomial 
system (3). In Step 2, we apply theorem \ref{main:height} to calculate the 
polynomial $h_{G_{(j)}}$. Since the $G_{(j)}$ all have an equal number 
of variables and equations (and none of the equations is of the form 
$0\!=\!0$), assertion (1) of theorem \ref{main:height} tells us that 
a particular $G_{(j)}$ has a complex root iff $h_{G_{(j)}}$ has positive 
degree. So it suffices to compute $h_{G_{(j)}}$ to check the feasibility of 
$G_{(j)}$. As for Step 4, note that thanks to theorem \ref{main:unired}, 
$G_{(j)}$ has a root in common with $F$ iff 
$\gcd\{h_{G_{(j)}},g_1(h_1,\dots,h_n),
\ldots,g_n(h_1,\ldots,h_n)\}$ 
has positive degree, where $h_1,\ldots,h_n$ are the polynomials 
corresponding to the application of theorem \ref{main:unired} to 
$G_{(j)}$. The preceding gcd and composition of univariate polynomials 
can be computed within $\cO(nk(n\log D)V_F\log^2 V_F)$ arithmetic 
operations via standard univariate polynomial algorithms \cite{binipan}, 
and we will soon see that this complexity is negligible compared to the 
work performed in the rest of our algorithm. 

Let us now check the correctness of our algorithm: 
Via lemma \ref{lemma:probe} and theorem \ref{main:height}, we see that Step 2 
gives a ``yes''  answer iff the intersection of $Z_{\twF}$ with a generic 
codimension $i$ flat is finite (and nonempty), where $\twF$ is an 
$n$-tuple of generic 
linear combinations of the $f_i$. Thus Step 2 gives a 
``yes'' answer iff $\dim Z_{\twF}\!=\!i$. 
Lemma \ref{lemma:gh} below tells us that $\dim Z_F\!=\!\dim Z_{\twF}$ if 
$\dim Z_F\!\geq\!1$. Otherwise, Step 5 correctly decides whether $Z_F$ is 
empty whenever $Z_F$ is finite. Thus the algorithm is correct. 

As for the complexity of our algorithm, letting $\cS$ (resp.\ $\cU$, $\cU'$) 
be the complexity of the corresponding application of lemma \ref{lemma:probe} 
(resp.\ theorems  \ref{main:height} and \ref{main:unired}), we immediately 
obtain a deterministic arithmetic complexity bound of 
\[ (n-2)\cS \ \ \mathrm{(All \ Executions \ of \ Step \ 1)} \] 
\[ +(n-2)(2k+1)\cU \ \ \mathrm{(All \ Executions \ of \ Step \ 2)} \] 
\[ +\cU'+\cO(n^2kV_F(\log^2 V_F)(\log D)) \ \ \mathrm{(Step \ 4)} \] 
(The complexity of the ``if'' statements in Steps 3 and 5 
is negligible.) Remark \ref{rem:probe} below tells us that 
$\cS\!=\!\cO((k+n^2)\log(m+n))$. Furthermore, in the proofs of theorems 
\ref{main:unired} and \ref{main:height} (cf.\ sections \ref{sub:unired} and 
\ref{sub:height}) later we will see that $\cU'\!=\cO(n\cU)$ and 
$\cU\!=\!\cO(V^3_F\cP_F)$. Since $k\!\geq\!m$, our overall complexity bound 
becomes $\cO(nk\cU+n\cS)\!=\!\cO(nkV^3_F\cP_F+n(k+n^2)\log(m+n))\!=
\cO(n^4kM^{2.376}_FV^5_F+n(k+n^2)\log(m+n))\!=\!\cO(
n^4kM^{2.376}_FV^5_F+nk\log(m+n))$. \qed 
\begin{rem} 
Note that if we somehow know that $\dim Z_F\!\geq\!1$, then we do not need 
assertion (2) of theorem \ref{main:height}, nor do we need to apply 
theorem \ref{main:unired}. We can thus pick any integral point (not equal to 
$\bO$) for $u_F$ and skip one of the steps of the proof of 
theorem \ref{main:height}. This removes a factor of $V^2_F$ from the 
first (usually dominant) summand of our complexity bound. \qed 
\end{rem} 

\begin{lemma}
\label{lemma:probe}  
Suppose $G(w,v)$ is a formula of the form 
\[\exists x_1\!\in\!\C \cdots \exists x_n\!\in\!\C \; 
(g_1(x,w,v)\!=\!0)\wedge \cdots \wedge (g_s(x,w,v)\!=\!0),\] 
where $g_1,\ldots,g_s\!\in\!\C[x_1,\ldots,x_n,w_1,\ldots,w_k, 
v_1,\ldots,v_r]$.  Then there is a sequence 
$v(1),\ldots,v(2k+1)\!\in\!\C^r$ such that 
for all $w\!\in\!\C^k$, the following statement holds: 
$G(w,v(j))$ is true for at least half of the $j\!\in\![2k+1]  
\Longleftrightarrow G(w,v)$ is true for a Zariski-open set of $v\!\in\!\C^r$. 
Furthermore, this sequence can be computed within $\log \sigma + (k+n+r)\log 
D$ arithmetic operations, where $\sigma$ (resp.\ $D$) is the maximum bit-size 
of any  coefficient (resp.\ maximum degree) of any $g_i$. \qed 
\end{lemma} 
\noindent 
The above lemma is actually just a special case of theorem 5.6 of 
\cite{koiran}. 
\begin{rem}
\label{rem:probe}
For the proof of theorem \ref{main:complex}, we have $s\!:=\!n$, 
$(g_1,\ldots,g_s)\!:=\!G_{(j)}$, $r\!:=\!(i+1)n\!\leq\!(n-1)n$, 
$v(j)\!=\!(\eps_1(j),\ldots,\eps_n(j),\eps_{(1,1)}(j),\ldots,
\eps_{(i,n)}(j))$, and we take $w$ to be the vector of coefficients 
of $F$. We thus obtain $\sigma\!=\!1$ and $D\!=\!m+i+1\!\leq\!m+n$. \qed 
\end{rem} 

\subsubsection{The Proof of Theorem \ref{main:height}}\mbox{}\\
\label{sub:height} 
\noindent
Curiously, precise estimates on coefficient growth in toric resultants 
are absent from the literature. So we supply such an estimate below. 
In what follows, we use $u_i$ in place of $u_{e_i}$. 
\begin{thm} 
\label{thm:growth} 
Following the notation of lemma \ref{lemma:respert}, suppose the coefficients 
of $F$ (resp.\ $F^*$) have absolute value bounded above by $c$ (resp.\ 
$c^*$) for all $i\!\in\![n]$ and $u_1,\ldots,u_n\!\in\!\C$. 
Also let $\|u\|\!:=\!\sqrt{u^2_1+\cdots+u^2_n}$ and let $\mu$ 
denote the maximal number of monomial terms in any $f_i$. 
Then the coefficient of $u^i_0$ in $\pert_{(F^*,\cA_{n+1})}$ has absolute 
value bounded above by \[ \frac{e^{13/12}}{\sqrt{\pi}}\sqrt{M_F+1}\cdot 
4^{M_F-i/2}\|u\|^{V_F-i}
(\sqrt{\mu}(c+c^*))^{M_F} 
\begin{pmatrix} V_F\\ i \end{pmatrix},\]
assuming that $\det \cM_{\bar{\cA}}\!\neq\!0$ under the substitution 
$(F-sF^*,u_0+u_1x_1+\cdots+u_nx_n) \mapsto \bar{F}$.  
(Note also that $\frac{e^{13/12}}{\sqrt{\pi}}\!\approx\!1.66691$.)   
\end{thm} 
\noindent 
{\bf Proof:} Let $c_{ij}$ denote the coefficient of 
$u^i_0s^j$ in $\det \cM_{\bar{\cA}}$, under the substitution
$(F-sF^*,u_0+u_1x_1+\cdots+u_nx_n) \mapsto \bar{F}$. 
Our proof will consist of computing an upper bound on $|c_{ij}|$, so 
we can conclude simply by maximizing over $j$ and then invoking a 
quantitative lemma on factoring. 

To do this, we first observe that one can always construct 
a toric resultant matrix with exactly $n_F$ rows corresponding 
to $f_{n+1}$ (where $\delta(Z_F)\!\leq\!n_F\!\leq\!V_F$), and the 
remaining rows corresponding to $f_1,\ldots,f_n$. 
(This follows from the algorithms we have already invoked 
in lemma \ref{lemma:respert}.) Enumerating how appropriate collections rows 
and columns can contain $i$ entries of $u_0$ (and $j$ entries 
involving $s$), it is easily verified that $c_{ij}$ is a sum of no more than 
$\begin{pmatrix} V_F \\ i \end{pmatrix} 
\begin{pmatrix} M_F-i \\ j \end{pmatrix}$ 
subdeterminants of $\cM_{\bar{A}}$ of size no greater than $M_F-i-j$. 
The coefficient $c_{ij}$ also receives similar contributions 
from some larger subdeterminants since the rows of $\cM_{\bar{\cA}}$ 
corresponding to $f_1,\ldots,f_n$ involve terms of the 
form $\eta+\nu s$. 

Via lemma \ref{lemma:multi} below, we can then derive 
an upper bound of 
\[\begin{pmatrix} V_F \\ i \end{pmatrix}\begin{pmatrix} M_F-i\\ j\end{pmatrix} 
\|u\|^{V_F-i}(\sqrt{\mu}(c+c^*))^{M_F-j} \] 
on $|c_{ij}|$. However, what we really need is an estimate on 
the coefficient $c_i$ of $u^i_0$ of $\pert_{(F^*,\cA_{n+1})}$, assuming 
the non-vanishing of $\det \cM_{\bar{\cA}}$. To 
estimate $c_i$, we simply apply lemma \ref{lemma:mignotte} below 
(observing that $\pert_{(F^*,\cA_{n+1})}$ is a divisor 
of an $M_F\times M_F$ determinant) to obtain an upper bound of 
\[\sqrt{M_F+1}\cdot 2^{M_F}\begin{pmatrix} V_F \\ i \end{pmatrix}\max_j\left\{
\begin{pmatrix} M_F-i\\ j\end{pmatrix}\right\} 
\|u\|^{V_F-i}(\sqrt{\mu}(c+c^*))^{M_F} \] 
on $|c_i|$. We can then finish via the elementary inequality $\begin{pmatrix} 
M_F-i \\ j \end{pmatrix}\!\leq\!\frac{e^{13/12}}{\sqrt{\pi}}2^{M_F-i}$,  
valid for all $j$ (which in turn is a simple corollary of Stirling's 
formula). \qed 

A simple result on the determinants of certain symbolic matrices, 
used above, is the following. 
\begin{lemma}
\label{lemma:multi} 
Suppose $A$ and $B$ are complex $N\times N$ matrices, where 
$B$ has at most $N'$ nonzero rows. Then the coefficient of 
$s^j$ in $\det(A+sB)$ has absolute value no greater than 
$\begin{pmatrix} N' \\ j \end{pmatrix}v^{N-j}(v+w)^j$, 
where $v$ (resp.\ $w$) is any upper bound on the 
Hermitian norms of the rows of $A$ (resp.\ $B$). \qed 
\end{lemma} 
\noindent 
The lemma follows easily by reducing to the case $j\!=\!0$, via 
the multilinearity of the determinant. The case $j\!=\!0$ is 
then nothing more than the classical {\bf Hadamard's lemma} 
\cite{mignotte}. 

The lemma on factorization we quoted above is the following. 
\begin{lemma} 
\cite{mignotte}
\label{lemma:mignotte}
Suppose $f\!\in\!\Z[x_1,\ldots,x_N]$ has total degree $D$ and
coefficients of absolute value $\leq\!c$. 
Then $g\!\in\Z[x_1,\ldots,x_N]$ divides $f \Longrightarrow$ the 
coefficients of $g$ have absolute value $\leq\!\sqrt{D+1}\cdot 2^Dc$. \qed
\end{lemma} 

We are now ready to prove theorem \ref{main:height}:\\
{\bf Proof of Theorem \ref{main:height}}:\\
By adjusting the number polynomials $m$ we can immediately assume 
that no $f_i$ is indentically zero. Furthermore, if $m\!=\!0$, 
we can clearly set $h\!:=\!0$. So we can also assume that $m\!\geq\!1$. 
We now consider three obvious cases. 

\noindent 
{\bf (The Case $\pmb{m\!=\!n}$):} 
The existence of an $h_F$ satisfying (0)--(5) will follow from 
setting $h_F(u_0)\!:=\!\pert_{(F^*,\cA_{n+1})}(u_0)$ for 
$\cA_{n+1}$ as in lemma \ref{lemma:respert}, $F^*$ as in lemma 
\ref{lemma:fill} below, and picking several $(u_1,\ldots,u_n)$ until a good 
one is found. Assertion (0) of theorem \ref{main:height} thus follows 
trivially. That the conclusion of lemma \ref{lemma:fill} implies assertion 
(1) is a consequence of \cite[Def.\ 2.2 and Main Theorem 2.1]{gcp}.  
 
To prove assertions (1)--(5) together we will then need to pick 
$(u_1,\ldots,u_n)$ subject to a final technical condition. In particular, 
consider the following method: 
Pick $\eps\!\in\![1+\begin{pmatrix} V_F \\ 2\end{pmatrix}]$ and set 
$u_i\!:=\!\eps^i$ for all $i\!\in\![n]$.
The worst that can happen is that  
a root of $h_F$ is the image two distinct points 
in $Z_F$ under the map $(\zeta_1,\ldots,\zeta_n) \mapsto 
u_1\zeta_1+\cdots+u_n\zeta_n$, thus obstructing assertion (2). (Whether this 
happens can easily be checked within $\cO(V_F\log V_F)$ arithmetic 
operations via a gcd calculation detailed in \cite[Sec.\ 5.2]{gcp}, 
after first finding the coefficients of $h_F$.) 
Otherwise, it easily follows from Main Theorems 2.1 and 2.4 of \cite{gcp}  
(and theorem \ref{main:unired} above and theorem \ref{thm:growth} below) 
that $h_F$ satisfies assertions (1)--(3) and (5). 

Since there are at most $\begin{pmatrix} 
V_F \\ 2\end{pmatrix}$ pairs of points $(\zeta_1,\zeta_2)$, 
picking $(u_1,\ldots,u_n)$ as specified above {\bf will} eventually 
give us a good $(u_1,\ldots,u_n)$. The overall arithmetic complexity of our 
search for $u_F$ and $h_F$ is, thanks to lemma \ref{lemma:respert},\\ 
$(\begin{pmatrix}V_F \\ 2 \end{pmatrix}+1) \cdot 
(V_F\cP_F+\cO(V_F\log V_F))$. This proves assertion (4), and we are done. 
\qed 
\begin{rem} 
Note that we never actually had to compute $V_F$ above: To pick a 
suitable $u$, we simply keep picking choices (in lexicographic order) with 
successively larger and larger coodinates until we find a suitable $u$. \qed 
\end{rem} 

\noindent 
{\bf (The Case $\pmb{m\!<\!n}$):} Take $f_{n+1}\!=\cdots =\!f_m\!=\!f_n$. 
Then we are back in the case $m\!=\!n$ and we are done. \qed 

\noindent 
{\bf (The Case $\pmb{m\!>\!n}$):} Here we employ an old trick: We substitute 
generic linear combinations of $f_1,\ldots,f_m$ for $f_1,\ldots,f_n$. 
In particular, set $\tf_i\!:=f_1+\eps_if_2+\cdots+\eps^{m-1}_if_m$ for 
all $i\!\in\![n]$. It then follows from lemma \ref{lemma:gh} below 
that, for generic $(\eps_1,\ldots,\eps_n)$, $Z_{\twF}$ is the union of $Z_F$ 
and a (possibly empty) finite set of 
points. So by the $m\!=\!n$ case, and taking into account the larger 
value for $c$ in our application of theorem \ref{thm:growth}, we are done. \qed 
\begin{rem} 
\label{rem:height}
Following the notation of theorem \ref{thm:growth}, we thus see that the 
asymptotic bound of assertion (3) can be replaced by an explicit bound of 
\[ \log\left\{\frac{e^{13/6}}{\pi}\sqrt{M_F+1}\cdot 
2^{V_F}4^{M_F}\left(\sqrt{n}\left(\begin{pmatrix}V_F \\ 2\end{pmatrix}
+1\right)^n\right)^{V_F} (c+1)^{M_F} \right\}\]
if $m\!\leq\!n$, or 
\[ \log\left\{\frac{e^{13/6}}{\pi}\sqrt{M_F+1}\cdot 
2^{V_F}4^{M_F}\left(\sqrt{n}\left(\begin{pmatrix}V_F \\ 2\end{pmatrix}
+1\right)^n\right)^{V_F} 
\sqrt{\mu}^{M_F}(m(mV_F+1)^{m-1}c+1)^{M_F}\right\} \]
for $m\!>\!n\!\geq\!1$. \qed 
\end{rem} 
\begin{lemma} 
\label{lemma:fill} 
Following the notation above 
let $\cA^*_i\!=\!\{\bO,e_1,\ldots,e_n\}\cup\bigcup^n_{j=1}\cA_j$ for all 
$i\!\in\![n]$ and $k^*\!:=\!n\#\cA_1$, where 
$\#$ denotes set cardinality. 
Also let $\cC^*$ be the coefficient vector of $F^*$. 
Then there is an $F^*$ such that (i) 
$\supp(F^*)\!\subseteq\!\cA^*$, (ii) 
$\cC^*\!=\!(1,\ldots,1)$, (iii) $F^*$ has exactly 
$V_F$ roots in $\Csn$ counting multiplicities, and 
(iv) $\det \cM_{\bar{\cA}}\!\neq\!0$ under the substitution 
$(F-sF^*,u_0+u_1x_1+\cdots+u_nx_n) \mapsto \bar{F}$. \qed 
\end{lemma} 
\noindent 
The above lemma is a paraphrase of \cite[Definition 2.3 and Main Theorem 
2.3]{gcp}. Furthermore, the deterministic arithmetic complexity of finding 
such an $F^*$ is dominated by $\cO(M_F\log n+n^2)$ \cite{new}, and can thus 
be ignored in our main bounds. 

\begin{lemma}
\label{lemma:gh}
Following the notation above, let $S\!\subset\!\C$ be any finite set 
of cardinality $\geq\!mV_F+1$. Then there is an 
$(\eps_1,\ldots,\eps_n)\!\in\!S^n$ such that every irreducible component of 
$Z_{\twF}$ is either an irreducible component of $Z_F$ or a point. \qed  
\end{lemma} 
\noindent 
The proof is essentially the same as the first theorem of \cite[Sec.\ 
3.4.1]{giustiheintz}, save that we use part (0) of theorem \ref{main:height} 
in place of B\'ezout's Theorem. 

\subsubsection{The Proof of Theorem \ref{main:size} }\mbox{}\\ 
Since we only care about the size of $x_i$, we can simply
pick $u_0\!=\!-1$, $u_i\!=1$, all other $u_j\!=\!0$, and apply
the polynomial $h_F$ from theorem \ref{main:height}.
(In particular, differing from the proof of theorem \ref{main:height},
we need not worry if our choice of $(u_1,\ldots,u_n)$ results in two distinct
$\zeta\!\in\!Z_F$ giving the same value for $\zeta_1u_1+\cdots+\zeta_nu_n$.)
Thus, by following almost the same proof as assertion (3) of theorem
\ref{main:height}, we can beat the height bound from theorem
\ref{main:height} by a summand of $\cO(n^2V_F\log D)$. \qed
\begin{rem}
\label{rem:size}
Via theorem \ref{thm:growth} (and a classic root
size estimate of Cauchy \cite{mignotte}), we easily see that the
asymptotic bound for $|\log|x_i||$ can be replaced by 
explicit quantities slightly better than those stated 
in remark \ref{rem:height}. In particular, it is clear from our 
last proof that we can simply replace  
the terms of the form $\sqrt{n}\left(\begin{pmatrix}V_F \\ 2\end{pmatrix}
+1\right)^n$ in the formulae from remark \ref{rem:height} 
by $\sqrt{2}$. \qed
\end{rem}

\subsubsection{The Proof of Theorem \ref{main:unired}}\mbox{}\\ 
\label{sub:unired} 
All portions, save assertion (8), follow immediately from 
\cite[Main Theorem 2.1]{gcp}. To prove assertion (8), we 
will briefly review the computation of $h_1,\ldots,h_n$ 
(which was already detailed at greater length in \cite{gcp}). 
Our height bound will then follow from some elementary  
polynomial and linear algebra bounds. 

In particular, recall the following algorithm for computing 
$h_1,\ldots,h_n$, given $h$ as in theorem \ref{main:height}: 
\begin{itemize} 
\item[{\bf Step 2}]{If $n\!=\!1$, set $h_1(\theta)\!:=\!\theta$ and stop.
Otherwise, for all $i\!\in\![n]$, let $q^-_i(t)$ be the square-free part of
$\pert_A(t,u_1,\ldots,u_{i-1},u_i-1, u_{i+1},\ldots,u_n)$.}
\item[{\bf Step 3}]{Define $q^\star_i(t)$ to be the square-free part of
$\pert_A(t,u_1,\ldots,u_{i-1},u_i+1,u_{i+1},\ldots,u_n)$ for all
$i\!\in\![n]$.}
\item[{\bf Step 4}]{For all $i\!\in\![n]$ and $j\!\in\!\{0,1\}$, let
$r_{i,j}(\theta)$ be the reduction of $\cR_j(q^-_i(t),
q^\star_i((\alpha+1)\theta-\alpha t))$ modulo $h(\theta)$. }
\item[{\bf Step 5}]{For all $i\!\in\![n]$, define
$g_i(\theta)$ to be the reduction of
$-\theta-\frac{r_{i,1}(\theta)}{r_{i,0}(\theta)}$ modulo $h(\theta)$. 
Then define $a_i$ to be the least positive integer so that 
$h_i(t)\!:=\!a_ig_i\!\in\!\Z[t]$. } 
\end{itemize}

Following the notation of the algorithm above, the polynomial 
$\cR_0(f,g)+\cR_1(f,g)t$ is known as the {\bf first subresultant} of $f$ and 
$g$ and can be computed as follows: Letting 
$f(t)\!=\!\alpha_0+\alpha_1t+\cdots+\alpha_{d_1}t^{d_1}$ and
$g(t)\!=\!\beta_0+\beta_1t+\cdots+\beta_{d_2}t^{d_2}$, consider the following
$(d_1+d_2-2)\times (d_1+d_2-1)$ matrix
\begin{small}
\[
\begin{bmatrix}
\beta_0 & \cdots & \beta_{d_2} & 0   & \cdots & 0 & 0 \\
0      & \beta_0 & \cdots & \beta_{d_2} & 0 & \cdots & 0\\
\vdots & \ddots & \ddots &  & \ddots & \ddots & \vdots \\
0      & \cdots & 0 & \beta_0 & \cdots & \beta_{d_2} & 0 \\
0      & 0  & \cdots & 0 & \beta_0 & \cdots & \beta_{d_2} \\
\alpha_0 & \cdots & \alpha_{d_1} & 0 & \cdots & 0 & 0 \\
0     & \alpha_0 & \cdots & \alpha_{d_1} & 0 & \cdots & 0\\
\vdots &  \ddots  & \ddots &    &  \ddots  & \ddots  & \vdots   \\
0      & \cdots & 0 & \alpha_0 & \cdots & \alpha_{d_1} & 0 \\
0      & 0  & \cdots & 0 & \alpha_0 & \cdots & \alpha_{d_1}
\end{bmatrix}
\]
\end{small}
\hspace{-\sh}with $d_1\!-\!1$ ``$\beta$ rows'' and $d_2\!-\!1$ ``$\alpha$
rows.'' Let $M^1_1$ (resp.\ $M^1_0$) be the submatrix obtained by
deleting the last (resp.\ second to last) column. We then define 
$\cR_i(f,g)\!:=\!\det(M^1_i)$ for $i\!\in\!\{0,1\}$.

Continuing our proof of Theorem \ref{main:unired}, we see that 
we need only bound the coefficient growth of the intermediate 
steps of our preceding algorithm. Thanks to theorem \ref{thm:growth}, 
this is straightforward: First note that 
$\sigma(q^-_i)\!=\!\log((V_F+1)\cdot 2^{V_F})+\sigma(\bar{h}_F)$, 
where $\bar{h}_F$ is the square-free part of $h_F$. (This follows trivially 
from expressing the coefficients of a univariate polynomial $f(t+1)$ 
in terms of the coefficients of $f(t)$.) Via lemma \ref{lemma:mignotte} 
we then see that $\sigma(\bar{h}_F)\!=\!\log(\sqrt{V_F+1}\cdot 2^{V_F})+
\sigma(h_F)$, and thus $\sigma(q^-_i)\!=\!\cO(\sigma(h_F))$. 
Similarly, $\sigma(q^\star_i)\!=\!\cO(\sigma(h_F))$ as well. 

To bound the coefficient growth when we compute $r_{i,j}$ note 
that the coefficient of $t_i$ in $q^\star_i(2\theta-t)$ 
is exactly $(-1)^i\sum^d_{j=i} \begin{pmatrix}j \\ i\end{pmatrix} 
(2\theta)^j\alpha_j$, where $\alpha_j$ is the coefficient of 
$t^j$ in $q^\star_i(t)$. Thus, via Hadamard's lemma again, 
we see that 
\[ |r_{i,j}(\theta)|\!\leq\!\left(\sqrt{V_F+1}\cdot 
e^{\sigma(h_F)}\right)^{V_F-1} \left(\sqrt{V_F+1}\cdot 
V_F2^{V_F}(2\theta)^{V_F}
e^{\sigma(h_F)}\right)^{V_F-1}\] for all $i,j$. Since $r_{i,j}$ is 
itself a polynomial in $\theta$ of degree $V_F(V_F-1)$, the 
last inequality then easily implies that 
$\sigma(r_{i,j})\!=\!\cO(V_F\sigma(h_F))$. 

To conclude, note that for any univariate polynomials $f,g\!\in\!\Z[t]$ 
with degree $\leq\!D$, $\sigma(fg)\!=\cO(\sigma(f)+\sigma(g)+\log D)$. 
Via long division it also easily follows that the 
quotient $q$ and remainder $r$ of $f/g$ satisfy $aq,ar\!\in\Z[t]$ 
and $\sigma(aq),\sigma(ar)\!=\!\cO(D(\sigma(f)+\sigma(g)))$, for some 
positive integer $a$ with $\log a\!=\!\cO(\sigma(g))$. 

So by assertion (3) of theorem \ref{main:height} we obtain 
$\log(a_i),\sigma(h_i)\!=\!\cO(V^2_F\sigma(h_F))$. \qed 
\begin{rem} 
\label{rem:denom}
An immediately consequence of our proof is that the 
asymptotic bound from assertion (8) can be replaced 
by the following explicit bound: 
\[ V_F\left\{(V_F-1)\left[\log\left(V_F(V_F+1)^{4} {64}^{V_F}\right)+
2\sigma(h_F)\right]+\sigma(h_F) \right\}+\sigma(h_F)+\log V_F. 
\text{ \ \qed} \] 
\end{rem}  

\subsubsection{The Proof of Theorem \ref{main:koi}}\mbox{}\\ 
\label{sub:koi}  
\noindent 
{\bf Proof of Part (a):} We first recall the following 
useful effective arithmetic Nullstellensatz of Krick, Pardo, and 
Sombra. 
\begin{thm}
\label{thm:cool} 
Suppose $f_1,\ldots,f_m\!\in\!\Z[x_1,\ldots,x_n]$ and 
$f_1\!=\cdots =\!f_m\!=\!0$ has {\bf no} roots in $\Cn$. 
Then there exist polynomials $g_1,\ldots,g_m\!\in\!\Z[x_1,\ldots,x_n]$ 
and a positive integer $a$ such that $g_1f_1+\cdots +g_mf_m\!=\!a$. 
Furthermore, \[ \log a\!\leq\!2(n+1)^3D V_F[\sigma(F)+\log m + 
2^{2n+4}D\log(D+1)]. \text{ \ \qed } \] 
\end{thm} 
\noindent 
The above theorem is a portion of corollary 3 from \cite{cool}. 

The proof of part (a) is then almost trivial: By assumption, 
theorem \ref{thm:cool} tells us that the mod $p$ reduction of $F$ 
has a root in $\Z/p\Z \Longrightarrow p$ divides $a$. Since 
the number of divisors of an integer $a$ is no more than 
$1+\log a$ (since any prime power other than $2$ is bounded below by 
$e$), we arrive at our desired asymptotic bound on $a_F$. \qed 
\begin{rem}
\label{rem:shebanga}
Following the notation of theorem \ref{main:koi}, 
we thus obtain the following explicit bound:  
\[ a_F\!\leq\!1+2(n+1)^3D V_F[\sigma(F)+\log m + 2^{2n+4}D\log(D+1)]. 
\text{ \ \qed } \] 
\end{rem}  

\noindent 
{\bf Proof of Part (b):} Recall the following version of the 
discriminant. 
\begin{dfn} 
\label{dfn:disc}
Given any polynomial
$f(x_1)\!=\!\alpha_0+\alpha_1x_1+\cdots+\alpha_Dx^D_1\!\in\!\Z[x_1]$
with all $|\alpha_i|$ bounded above by some integer $c$, define the
{\bf discriminant of} $\mathbf{f}$, $\pmb{\Delta_f}$, to be
$\frac{(-1)^{D(D-1)/2}}{\alpha_D}$ times the following
$(2D-1)\times (2D-1)$ determinant:
\begin{small}
\[ \det
\begin{bmatrix}
\alpha_0 & \cdots & \alpha_D & 0 & \cdots & 0 & 0 \\
0     & \alpha_0 & \cdots & \alpha_D & 0 & \cdots & 0\\
\vdots &  \ddots  & \ddots &    &  \ddots  & \ddots  & \vdots   \\
0      & \cdots & 0 & \alpha_0 & \cdots & \alpha_D & 0 \\
0      & 0  & \cdots & 0 & \alpha_0 & \cdots & \alpha_D \\
\alpha_1 & \cdots & D\alpha_D & 0 & \cdots & 0 & 0 \\
0     & \alpha_1 & \cdots & D\alpha_D & 0 & \cdots & 0\\
\vdots &  \ddots  & \ddots &    &  \ddots  & \ddots  & \vdots   \\
0      & \cdots & 0 & \alpha_1 & \cdots & D\alpha_D & 0 \\
0      & 0  & \cdots & 0 & \alpha_1 & \cdots & D\alpha_D 
\end{bmatrix},
\]
\end{small} 
\noindent
\mbox{}\hspace{-.15cm}where the first $D-1$ (resp.\ last $D$) rows correspond 
to the coefficients of $f$ (resp.\ the derivative of $f$). \qed
\end{dfn} 

Our proof of part (b) begins with the following observation. 
\begin{thm} 
\label{thm:oyster} 
Following the notation of section \ref{sec:rat}, 
suppose $f\!\in\!\Z[x_1]$ is a square-free polynomial of 
degree $D$ with exactly $i_f$ factors over $\Q[x_1]$. 
Then the truth of GRH implies that  
\[ \left|i_f\pi(t)-N_f(t)\right|\!<\!2\sqrt{t}(D\log t+\log 
|\Delta_f|) +D\log |\Delta_f|,  \] 
for all $t\!>\!2$. \qed  
\end{thm} 
\noindent 
A slightly less explicit version of the above theorem appeared 
in \cite[Thm.\ 9]{hnam}, and the proof is almost the same as that of an  
earlier result of Adleman and Odlyzko for the case $i_f\!=\!1$ 
\cite[Lemma 3]{amo}. (See also \cite{weinberger}.) The only new 
ingredient is an explicit version of the effective Chebotarev density theorem 
due to Oesterl\'e \cite{oyster}. (Earlier versions of theorem \ref{thm:oyster} 
did not state the asymptotic constants explicitly.) 

The proof of part (b) is then essentially a chain of elementary 
analytic bounds which flows from applying theorem \ref{thm:oyster} 
to the polynomial $h_F$ from theorem \ref{main:complex}. However, a 
technicality which must be considered is that $h_F$ might not be 
square-free (i.e., $\Delta_{h_F}$ may vanish). This is easily taken care of by 
an application of the following immediate corollary of lemmata  
\ref{lemma:multi} and \ref{lemma:mignotte}. 
\begin{cor} 
\label{cor:disc}
Following the notation above, let $g$ be the square-free part 
of $f$ and let $D'$ be the degree of $g$. Then 
$\log |\Delta_g|\!\leq\!D'(D\log 2+\log(D'+1)+\log c)$. \qed 
\end{cor}

Another technical lemma we will need regards the existence of primes 
interleaving a simple sequence. 
\begin{lemma} 
\label{lemma:pain} 
The number of primes in the open interval $(At^3,A(t+1)^3)$ 
is at least $\lfloor \frac{1}{12}\cdot\frac{At^2} {\log t+\log A}\rfloor$, 
provided 
$A,t\!>\!e^5\!\approx\!148.413$. \qed 
\end{lemma} 
\noindent 
This lemma follows routinely (albeit a bit tediously) from theorem 8.8.4 of 
\cite{bs}, which states that for all $t\!>\!5$, the $t^\thth$ prime lies in 
the open interval $(t\log t,t(\log t+\log\log t))$. 

The key to proving theorem \ref{main:koi} is then to find small constants 
$t_0$ and $A_F$ such that $N_F(A_F(t+1)^3-1)-N_F(A_Ft^3)\!>\!1$ for all 
$t\!\geq\!t_0$. 

Via theorems \ref{main:height} and \ref{main:unired}, and a consideration of 
the primes dividing the $a_i$ (the denominators in our rational 
univariate representation of $Z_F$), it immediately follows that 
$|N_F(t)-N_{h_F}(t)|\!\leq\!V_F\sum^n_{i=1}(\log a_i+1)$, for all $t\!>\!0$. 
We are now ready to derive an inequality whose truth will 
imply $N_F(A_F(t+1)^3-1)-N_F(A_Ft^3)\!>\!1$: 
By theorem \ref{thm:oyster}, lemma \ref{lemma:pain}, the triangle inequality, 
and some elementary 
estimates on $\log t$, $t^3$, and their derivatives, it suffices to 
require that $A_Ft^2$ strictly exceed $12(\log A_F+\log t)$ 
times the following quantity: 
\[2(1+\sqrt{2})\sqrt{3A_Ft^3}[V_F(\log(3A_Ft^3)+1)+\log |\Delta_g|]+
V_F\left(\log|\Delta_g| + \sum^n_{i=1}\log a_i +n\right)+1, \] 
for all $t\!>\max\{t_0,e^5\}$, where $g$ denotes the square-free part of 
$h_F$. (Note that we also used the fact that $i_g\!\geq\!1$.) 

A routine but tedious estimation then shows that 
we can actually take $t_0\!=\!1296(\frac{1+\log 3}{3}+\log 1296)\!
\approx\!4963040.506$, and $A_F$ as in the statement of part (b). 
Careful accounting of the estimates then easily yields the explicit 
upper bound for $A_F$ we state below. \qed 
\begin{rem} 
\label{rem:shebangb}
The constant $1296(\frac{1+\log 3}{3}+\log 1296)$ arises from trying to find 
the least $t$ for which $t^2\!\geq\!\alpha \log^4t$, where, roughly 
speaking, $\alpha$ ranges over the constants listed in the 
expressions for $A_F,B_F,C_F,D_F$ below. 
\[ A_F\!\leq\!\lceil 1296B^2_F\log^4B_F+36C^2_F\log^2C_F+2D_F\log D_F\rceil, \]
where 
\[ B_F\!:=\!72\sqrt{3}(1+\sqrt{2}) V_F, \ 
C_F\!:=\!24\sqrt{3}(1+\sqrt{2})\log|\Delta_g| +2, \text{ \ and} \]  
\[ D_F\!:=\!12V_F\left(\log|\Delta_g|+\sum^n_{i=1}\log a_i+n\right)+13. 
\text{ \ \qed } \] 
%
\end{rem} 

\subsection{Proofs of Our Results Over $\Z$: Theorems \ref{main:equiv}, 
\ref{main:25}, and \ref{main:h3}}\mbox{}\\
\label{sub:proofint} 
The proof of theorems \ref{main:equiv} and \ref{main:25} rely on a refined 
version of Siegel's theorem (theorem \ref{thm:siegel} stated earlier in 
section \ref{sec:int}) and an algorithmic result on factoring polynomials 
over $\C$ (lemma \ref{lemma:fac} below). 
The proof of theorem \ref{main:h3} will mainly use the tools we developed for 
our results over $\C$ from section \ref{sec:complex}, and is a streamlined 
version of the proof from \cite{tcs}.

\subsubsection{The Proof of Theorem \ref{main:equiv}}\mbox{}\\ 
\noindent
{\bf ($\mathbf{\Longrightarrow}$):} Simply apply whatever 
algorithm one has for $\biggy_\N$ to the polynomial  
$f(-x,-y)f(-x,y)f(x,-y)f(x,y)$ to obtain the value of  
$\biggy_\Z(f)$. \qed 

\noindent 
{\bf ($\mathbf{\Longleftarrow}$):} First calculate 
$b\!:=\!\biggy_\Z(f)$. If $b\!<\!\infty$ then we can 
simply enumerate {\bf positive} integral points until we 
at last know $\biggy_\N(f)$. (This can of course be mind-bogglingly slow, 
but is nevertheless a Turing-machine algorithm which is guaranteed to 
terminate.) 

If $b\!=\!\infty$ then let us do the following: Replace $f$ by 
its square-free part. (This can be done within $\nc$ via, say,  
lemma \ref{lemma:fac} below.) 
Then note that any irreducible component of $Z_f$ containing infinitely many 
integral points must be defined over $\Z$. (Otherwise, the action of 
$\mathrm{Gal}(\bar{\Q}/\Q)$ would imply that every integral point has 
multiplicity $>\!1$ --- a contradiction, since the number of 
singular points of a curve is always finite.) So we may also 
assume that $Z_f$ is geometrically irreducible. (Indeed, 
we can find all the irreducible components of $Z_f$ within 
$\nc$ via lemma \ref{lemma:fac}.) 

Theorem \ref{thm:silvb} then tells us that $\biggy_\N(f)\!=\!\infty 
\Longleftrightarrow Z_f$ has unbounded intersection with the 
the (open) first quadrant. To decide the latter question, one first finds the 
largest real critical value of the projection $(x,y)\mapsto x+y$, restricted 
to the intersection of $Z_f$ with the first quadrant. (Since we are 
restricting to the first quadrant, one must also consider the 
image of the intersection of $Z_f$ with the coordinate axes under 
this projection as well.) This reduces to finding 
the $(\zeta_1,\zeta_2)$ which maximizes $\zeta_1+\zeta_2$, 
where $(\zeta_1,\zeta_2)$ is either a positive real roots of the polynomial 
system $(f,\frac{\partial f}{\partial x} +\frac{\partial f}{\partial y})$, 
or a point in $Z_f\cap\{xy\!=\!0\}$. 
Thanks to theorems \ref{main:height} and 
\ref{main:unired}, and a fast root approximation algorithm from 
\cite{neffreif}, this can be done within $\nc$. 

To conclude, if there is no critical value, we simply check (via the 
techniques just mentioned) if the polynomial system $(f,x+y-1)$ 
has a positive real root. It is then easily checked that 
this system has a root iff $Z_f$ has unbounded intersection 
with the first quadrant. Otherwise, one performs the same 
check with the polynomial system $(f,x+y-\zeta_1-\zeta_2-1)$ instead. 
So we are done. \qed 

\subsubsection{The Proof of Theorem \ref{main:25}}\mbox{}\\ 
First note that as in our last proof, we can use 
lemma \ref{lemma:fac} to reduce (within $\nc$, relative to the 
dense encoding) to the case where $Z_f$ is geometrically 
irreducible.  

Our algorithm then proceeds as follows: Compute the genus 
of $Z_f$. (By \cite{ks}, this can actually be done within $\nc$ as well.) 
If the genus is positive then theorem \ref{thm:siegel} tells us that 
there are only finitely many integral points and we are done. 
Similarly, via \cite{neffreif}, condition (c) of theorem \ref{thm:siegel} 
can be checked within $\nc$. 

So we may now assume that $Z_f$ satisfies condition (c) and has genus zero. 
Find all {\bf positive} integral singular points of $Z_f$. 
(By theorems \ref{main:height}, \ref{main:unired}, and \ref{thm:lenstra}, 
this can also be done within $\nc$.) Call 
these points $\{(\alpha_1,\beta_1),\ldots,(\alpha_N,\beta_N)\}$. 
Then form the polynomial $g(x,y,t)\!:=\!
(x-\alpha_1)^2+(y-\beta_1)^2+\cdots+(x-\alpha_N)^2+(y-\beta_N)^2-t$.
Clearly, $Z_f$ has a nonsingular integral point iff the curve 
$Z_{(f,g)}\!\subset\!\C^3$ has a positive integral point. 
Furthermore, since $Z_f$ has a rational parametrization, the curve 
$Z_{(f,g)}$ admits one as well. Thus $Z_{(f,g)}$ is irreducible and has genus 
zero too.   

So assuming $\ratcurve(3)$ is decidable, theorem \ref{thm:siegel} 
tells us that we can decide whether $Z_f$ has infinitely many 
integral points. Converting this to the decidability of 
$\htp^\infty(2)$ is a simple matter, thanks to theorem \ref{thm:silvb} 
and an application of theorem \ref{main:unired} already detailed 
in our last proof. \qed  

\begin{lemma} 
\label{lemma:fac} 
\cite{bcgw} 
Suppose $f\!\in\!\Q[x_1,\ldots,x_n]$ and $n$ is a constant.  
Then, relative to the dense encoding, we can find all factors of 
$f$ over $\C[x_1,\ldots,x_n]$ within $\nc$. Furthermore, 
every factor is given as a polynomial in $\Q[\alpha][x_1,\ldots,x_n]$, where 
the minimal polynomial of $\alpha$ is also part of the output. \qed 
\end{lemma} 

\subsubsection{The Proof of Theorem \ref{main:h3}}\mbox{}\\ 
It suffices to show that the truth of both conditions 
implies the existence of an algorithm for $\exists\exists\forall\exists$ 
(with all quantifiers ranging over $\N$), thus contradicting 
the aforementioned result of Matiyasevich and Robinson. 

So assuming the truth of (1) and (2), let us construct such an 
algorithm. First note the following fact. 
\label{sub:h3}
\begin{lemma} 
\label{lemma:bound}
Following the notation above, 
let \[\Sigma_f\!:=\!\{(u_0,v_0)\!\in\!\C^2 \; \;  
| \; \; \{(x,y)\!\in\!\C^2 \; | 
\; f(u_0,v_0,x,y)\!=\!0\} 
\text{ \ has \ a \ genus \ zero \ component}\}.\] 
Also let $\Xi_f$ denote the set of  
$(u_0,v_0)\!\in\!\N^2$ such that $\forall x \; \exists y \; 
f(u_0,v_0,x,y)\!=\!0$. 
Then $\Xi_f\!\subseteq\!\Sigma_f\cap\Z^2$, whether all quantifiers range over 
$\N$ or $\Z$. 
\end{lemma} 
\noindent 
{\bf Proof of the Lemma:} 
By theorem \ref{thm:siegel}, 
$\forall x\;  \exists y\;  f(u_0,v_0,x,y)\!=\!0 
\Longrightarrow Z_f\cap\{(u,v)\!=\!(u_0,v_0)\}$ contains a curve of 
genus zero (whether the quantification is over $\N$ or $\Z$). 
So we are done. \qed 

Continuing the proof of theorem \ref{main:h3}, consider the following 
algorithm for $\exists\exists\forall\exists$: First decide whether 
$Z_f$ contains a specially ruled surface. (That this is Turing-decidable 
was already observed in \cite{tcs}.) If so, simply apply any 
algorithm for statement (2) to decide the prefix 
$\exists\exists\forall\exists$. 

Otherwise, $\Sigma_f$ is the (possibly empty) union of a finite point set and 
a collection of curves of positive genus. Via algorithms already observed in 
\cite{tcs}, the defining polynomials for all these points and curves are 
Turing-computable. So via theorem \ref{main:unired}, and statement 
(1), the worst we need do is enumerate integral points on 
several curves of positive genus. So although our algorithm may be 
very slow, we have succeeded in deriving a contradiction, and we are done. \qed 

\begin{rem} 
The usual definition of genericity in computational algebra is 
stronger than the one we gave earlier: 
A statement involving a set of parameters $\{c_1,\ldots,c_N\}$ 
holds {\bf generically} iff the statement is true for all 
$(c_1,\ldots,c_N)\!\in\!\C^N$ outside of some {\bf a priori fixed} 
algebraic hypersurface. That this version of genericity implies 
the simplified version mentioned earlier in our theorems is 
immediate from Schwartz' Lemma \cite{schwartz}. Any 
statement claimed to be true generically in this paper still 
holds under this stronger notion. \qed 
\end{rem} 

\section{Acknowledgements} 
The author would like to express his deep gratitude to the 
organizers of this conference for their generous invitation. 
The author also thanks Felipe Cucker, Ioannis Emiris, Teresa Krick, 
Francois Loeser, Gregorio Malajovich, Luis-Miguel Pardo-Vasallo, 
Steve Smale, and Martin Sombra  
for some very useful discussions, in person and via e-mail. Many of the 
results presented in this paper would have been weaker, were it not 
for the wonderful atmosphere of the Hilbert 10 conference in Gent. 

I dedicate this paper to Steve Smale. 

\section*{Appendix: How the Examples Were Computed}  
\label{sec:app} 
Here we reveal some further details on the computations underlying 
our examples. All of the computations in this paper were performed 
on a Sun 4u Computeserver, named Kronecker, at MIT. The version 
of {\tt Maple} used was {\tt Maple V Release 5}. 

The univariate reduction, $P(u)$, for our first $3\times 3$ polynomial 
system is a nonzero constant multiple of the sparse 
resultant of $f_1$, $f_2$, $f_3$, and $u-xyz$. The following {\tt 
Maple} code is how the computation was performed: 
\begin{quote} 
{\tt 
\noindent
with(linalg);\\

\noindent
f:=144+2*x-3*y\^{}2+x\^{}7*y\^{}8*z\^{}9;\\
g:=-51+5*x\^{}2-27*z+x\^{}9*y\^{}7*z\^{}8;\\
h:=7-6*x+8*x\^{}8*y\^{}9*z\^{}7-12*x\^{}8*y\^{}8*z\^{}7;\\
k:=u-x*y*z;\\

\noindent
r1:=factor(resultant(f,k,x)):\\
r2:=factor(resultant(g,k,x)):\\
r3:=factor(resultant(h,k,x)):\\

\noindent 
rr1:=op(4,r1):\\
rr2:=op(4,r2):\\
rr3:=op(3,r3):\\

\noindent
s1:=factor(resultant(rr1,rr3,z)):\\
s2:=factor(resultant(rr2,rr3,z)):\\

\noindent
ss1:=op(4,s1):\\
ss2:=op(3,s2):\\

\noindent
t:=factor(resultant(ss1,ss2,y)):\\
univar:=op(3,t);\\
} 
\end{quote} 

We also note that our choice for $P(u)$ was a bit sneaky: 
instead of finding a polynomial whose roots were linear 
projection of the roots of $F$, we found a polynomial whose 
roots were a {\bf monomial map} of the roots of $F$. This 
additional flexibility is useful in practice, and it is also 
possible to improve our quantitative results along these lines. 
These improvements will be detailed in later work, and we 
also point out that other applications of such nonlinear projections 
have appeared in earlier work of the author \cite{esa}. 

As for the mixed volume calculation, we used a {\tt C} implementation 
by Ioannis Emiris (publically available at\\ 
{\tt http://www.inria.fr/saga/logiciels/emiris/soft\_geo.html}).
That the mixed volume equals the number of roots in $\C^3$ 
follows easily from the fact that all the polynomials have a 
nonzero constant term, and an exactness condition for Bernshtein's 
Theorem (see, e.g., \cite{bernie} or \cite[Main Theorem 2]{jpaa}). 
Verifying the latter condition amounts to checking whether a 
product of toric resultants vanishes and for the sake of brevity 
we omit this calculation. In any case, it is easily checked that 
$M_F\leq e^{3+\frac{1}{8}}\cdot\frac{243}{\sqrt{4}}+(3\cdot 9+2)^3-(3\cdot 
9+1)^3
\approx 5202.327253$ for our example, via lemma \ref{lemma:respert}. 
(In practice, the true value of $M_F$ is 
typically {\bf much} smaller than the upper bound from 
lemma \ref{lemma:respert}.)   

By a stroke of luck, the polynomial $P$ is irreducible over 
$\Q$, so we immediately obtain that $F$ has exactly $145$ 
{\bf distinct} complex roots. Furthermore, we obtain that for 
any subfield $K\!\subseteq\!\C$, every root of $P$ in $K$ 
is the image of a unique root of $F$ in $K^3$. So we also obtain that $F$ has 
no rational roots. Via the {\tt realroot} command of {\tt Maple} 
(which employs {\bf Sturm sequences} \cite{marie}), we similarly 
obtain the number of real roots of $F$. 

As for the comparison with Gr\"obner bases, we 
simply invoked the following {\tt Maple} commands: 

\begin{quote} 
{\tt 
\noindent
f:=144+2*x-3*y\^{}2+x\^{}7*y\^{}8*z\^{}9;\\
g:=-51+5*x\^{}2-27*z+x\^{}9*y\^{}7*z\^{}8;\\
h:=7-6*x+8*x\^{}8*y\^{}9*z\^{}7-12*x\^{}8*y\^{}8*z\^{}7;\\
k:=u-x*y*z;\\

\noindent 
with(Groebner);\\ 
univpoly(u,[f,g,h,k]);\\ 
}
\end{quote} 

The larger time bound given was actually the amount of time 
{\tt Maple} spent calculating a univariate reduction via 
Gr\"obner bases, until the author's remote connection to 
{\tt Kronecker} was terminated. 

\footnotesize
\bibliographystyle{acm}

\end{document}